\documentclass[11pt]{amsart}
\usepackage{mathrsfs,amsthm,amssymb,amsmath,url}

\theoremstyle{plain}
    \newtheorem{thm}{Theorem}
    
    \newtheorem{lem}[thm]{Lemma}
    \newtheorem{prop}[thm]{Proposition}
    \newtheorem{cor}[thm]{Corollary}

\theoremstyle{definition}
    \newtheorem{defn}[thm]{Definition}
    
\theoremstyle{remark}
    \newtheorem{rem}[thm]{Remark}

\newcommand{\nin}{\notin}

\DeclareMathOperator{\Const}{Const}
\newcommand{\rest}{\upharpoonright}

\newcommand{\sm}{{\setminus}}

\newcommand{\cmu}{($\mu$) }
\newcommand{\cnu}{($\nu$) }

\newcommand{\csig}{($\sigma$) }

\newcommand{\clamp}{($\lambda$') }
\newcommand{\cchi}{($\chi$) }

\newcommand{\cscont}{(scont) }

\newcommand{\cl}[1]{\langle #1 \rangle}
 \DeclareMathOperator{\supp}{supp}

\renewcommand{\O}{{\mathscr O}}
\newcommand{\On}{{\mathscr O}^{(n)}}

\newcommand{\Oo}{{\mathscr O}^{(1)}}

 \DeclareMathOperator{\pol}{Pol}

\newcommand{\gwithS}{\cl{\{g\}\cup\S}}
\newcommand{\fwithS}{\cl{\{f\}\cup\S}}
\newcommand{\inv}{^{-1}}

\newcommand{\C}{{\mathscr C}}

\newcommand{\F}{{\mathscr F}}

\newcommand{\I}{{\mathscr I}}
\newcommand{\A}{{\mathscr A}}

\newcommand{\G}{{\mathscr G}}

\newcommand{\K}{{\mathscr K}}
\renewcommand{\H}{{\mathscr H}}
\renewcommand{\c}{X\sm }
\newcommand{\sub}{\subseteq}

\renewcommand{\S}{{\mathscr S}}

\newcommand{\uo}{^{(1)}}
\newcommand{\ut}{^{(2)}}

\newcommand{\To}{\rightarrow}

\newcommand{\ef}{\varepsilon_f}
\newcommand{\efp}{\varepsilon_f'}

\author[H.\,Machida]{Hajime Machida}
\address{Department of Mathematics\\Hitotsubashi University\\Naka 2-1\\Kunitachi, Tokyo 186-8601, Japan}
\email{machida@math.hit-u.ac.jp}

\author[M.\,Pinsker]{Michael Pinsker}
\address{Algebra\\TU Wien\\Wiedner Hauptstra\ss e 8-10/104\\A-1040 Wien, Austria}
\email{marula@gmx.at} \urladdr{http://www.dmg.tuwien.ac.at}
\thanks{The second author is grateful for support through the Postdoctoral Fellowship of the
Japan Society for the Promotion of Science (JSPS)}
\title[Minimal clones above the permutations]
    {
        The minimal clones above the permutations
    }

\subjclass{Primary 08A40;secondary 08A05}

\keywords{clone lattice, clones containing the permutations,
monoids, minimal clones, intervals of the clone lattice}

\begin{document}

    \begin{abstract}
        We determine the atoms of the interval of the clone lattice consisting of those clones which
        contain all permutations, on an infinite base set. This is
        equivalent to the description of the atoms of the lattice
        of transformation monoids above the permutations.
    \end{abstract}

   \maketitle


    \section{The problem and the result}
        Let $X$ be an infinite set of cardinality $\kappa=\aleph_\alpha$,
        let $\O$ be the set of all finitary
        operations on $X$, and for all natural numbers $n\geq 1$ let
        $\On$ be the set of $n$-ary operations
        on $X$. A set of operations $\C\sub\O$ is called a
        \emph{clone} if and only if it is closed under composition
        of functions and contains all projections, i.e. the
        functions $\pi^n_k\in\On$ satisfying
        $\pi^n_k(x_1,\ldots,x_n)=x_k$, for all $n\geq 1$ and all
        $1\leq k\leq n$. Ordering the set of all clones on $X$ by
        set-theoretical inclusion, one obtains a complete
        algebraic lattice $Cl(X)$. The cardinality of $Cl(X)$ is
        easily seen to equal
        $2^{2^\kappa}$, and the lattice seems to be too complicated
        to ever be fully described.

        Therefore, it has been tried to investigate
        interesting parts of $Cl(X)$, such as the atoms, referred to as  \emph{minimal} clones,
        or the dual atoms, called \emph{maximal} or
        \emph{precomplete} clones. However, at least on infinite
        $X$, even describing the minimal clones or the maximal
        clones seems unrealistic, since despite considerable efforts the minimal
        clones are not even known
        in the much smaller clone lattice over a finite base set,
        and since the number of maximal clones on an infinite base
        set has been proven to equal $2^{2^\kappa}$
        (\cite{Ros76}, see also \cite{GS02}). Successful research
        has been done on intervals of the clone lattice, for
        example on $[\cl{\Oo},\O]$ in \cite{Gav65}, \cite{GS02}, \cite{GSveryMany}, and
        \cite{Pin032}, where $\cl{\Oo}$ denotes the clone generated by $\Oo$, and on the interval above
        the clone of idempotent functions in \cite{GSlargeIntervals}. Several results could also be obtained on the
        intervals $[\cl{\S},\O]$ and $[\cl{\S},\cl{\Oo}]$, where $\S$ is the
        monoid
        of all permutations on $X$: In \cite{Hei02} a complete list of the maximal
        clones of $[\cl{\S},\O]$, and in \cite{Gav65} one of
        the clones maximal in $[\cl{\S},\cl{\Oo}]$ were provided on a
        countably infinite base set $X$
        (the latter one being a list of monoids, since clones
        below $\cl{\Oo}$ consist of functions which depend on
        at most one variable and therefore correspond
        to monoids in an obvious way; we shall for this reason drop the brackets and talk about
        the interval $[\S,\Oo]$ of the monoid lattice).
        In \cite{Pin033}, the author extended the first result
        to sets of all regular cardinalities, and the second
        result to all infinite sets. It turned out that there
        exist $\max\{|\alpha|,\aleph_0\}$ maximal clones in
        $[\cl{\S},\O]$, and $2\cdot|\alpha|+5$ maximal monoids in
        $[\S,\Oo]$. Those numbers are relatively small considering the size of the clone lattice
        (or the monoid lattice, which is as large as the clone lattice), but the author proved in \cite{Pin041}
        that the cardinality of $[\S,\Oo]$ is
        $2^{2^{\max\{|\alpha|,\aleph_0\}}}$, so rather large.

        In this article, we determine all clones minimal in
        $[\cl{\S},\O]$. It turns out quickly that all such clones are in fact
        monoids, that is, they only contain functions depending on at
        most one variable. Therefore, the problem reduces to
        finding the minimal monoids of $[\S,\Oo]$, which is
        interesting in itself. We will see that there exist $\max\{|\alpha|,\aleph_0\}$ such
        monoids. Surprisingly, this implies that if $|X|<\aleph_\omega$, in particular on countably
        infinite $X$, there exist only finitely many maximal but infinitely many minimal elements in $[\S,\Oo]$.

        For a monoid $\G\sub\Oo$, define $\pol(\G)$ to consist of
        all $f\in\O$ for which $f(g_1,\ldots,g_n)\in\G$ whenever
        $g_1,\ldots,g_n\in\G$. Call a clone $\C$ \emph{collapsing} iff it is uniquely
        determined by its unary part $\C\cap\Oo$, that is, there exist no
        other clones with the same unary part. Equivalently, $\C$
        is collapsing iff all functions in $\pol(\C\uo)$
        are \emph{essentially unary},
        that is, they depend on at most one variable.
        \begin{lem}
            $\cl{\S}$ is collapsing.
        \end{lem}
        \begin{proof}
            Let $f\in\pol(\S)\cap \O\ut$. Then $\gamma(x)=f(x,x)$ is a permutation.
            Now let $a,b\in X$ be distinct. There exists $c\in X$ with
            $\gamma(c)=f(a,b)$. If $c\nin\{a,b\}$, then we can
            find $\alpha,\beta\in\S$ with $\alpha(a)=a$,
            $\alpha(b)=c$, $\beta(a)=b$, and $\beta(b)=c$. But
            then
            $f(\alpha,\beta)(a)=f(a,b)=f(c,c)=f(\alpha,\beta)(b)$,
            so $f(\alpha,\beta)(x)$ is not a permutation. Thus,
            $c\in\{a,b\}$, so we have shown that
            $f(x,y)\in\{f(x,x),f(y,y)\}$ for all $x,y\in X$.\\
            Next we claim that for all $a,b\in X$, if
            $f(a,b)=f(a,a)$, then $f(b,a)=f(b,b)$. Indeed, consider the permutation $\alpha$
            which has a cycle $(ab)$. Then $f(a,\alpha(a))=f(a,b)=f(a,a)$,
            so $f(b,\alpha(b))=f(b,a)$ has to be different from $f(a,a)$, because otherwise the function $f(x,\alpha(x))$ is not
            injective. Therefore, $f(b,a)=f(b,b)$.\\
            Assume without loss that $f(a,b)=f(a,a)$, for some distinct $a,b\in X$.
            We first claim that $f(a,c)=f(a,a)$ for all $c\in X$.
            For assume not; then $f(a,c)=f(c,c)$, and therefore $f(c,a)=f(a,a)$.
            Let $\beta\in\S$ map $a$ to $b$ and $c$ to $a$.
            Then $f(a,\beta(a))=f(a,b)=f(a,a)$, but also $f(c,\beta(c))=f(c,a)=f(a,a)$, a contradiction since $f$
            preserves $\S$. Hence, $f(a,c)=f(a,a)$ for all $c\in
            X$.\\
            Now if $f(\tilde{a},\tilde{b})\neq f(\tilde{a},\tilde{a})$ for some $\tilde{a},\tilde{b}\in X$,
            then the conditions $f(\tilde{a},\tilde{b})=f(\tilde{b},\tilde{b})$ but
            $f(a,\tilde{b})=f(a,a)\neq f(\tilde{b},\tilde{b})$ lead to a similar
            contradiction. Hence, $f(x,y)=f(x,x)$ for all $x,y\in
            X$, and we have shown that $f$ depends on at most one variable. Since $f\in\pol(\S)\cap\O\ut$ was arbitrary,
            all binary functions of $\pol(\S)$
            are essentially unary. By a result of Grabowski
            \cite{Gra97}, this implies that $\cl{\S}$ is collapsing. (The
            mentioned result was proved for finite base sets but
            the same proof works on infinite sets.)
        \end{proof}
        \begin{lem}
            If $\C$ is a clone that is minimal in $[\cl{\S},\O]$, then it contains only
            essentially unary functions.
        \end{lem}
        \begin{proof}
            Since $\C\supsetneq\cl{\S}$ we have
            $\C\uo\supseteq\S$. If $\C\uo=\S$, then $\C=\cl{\S}$
            since $\cl{\S}$ is collapsing, contradicting the
            assumption that $\cl{\S}$ be a proper subset of $\C$.
            Therefore, $\C\uo$ properly contains $\S$, so
            $\C=\cl{\C\uo}$ since $\C$ is minimal above $\S$.
            Hence $\C$ contains only essentially unary functions.
        \end{proof}

        By the preceding lemma, when looking for the minimal clones above $\cl{\S}$, it suffices to determine the
        minimal monoids above $\S$. Clearly, such monoids are
        generated by a single non-permutation together with $\S$.
        We call functions which generate a minimal monoid above the permutations $\S$-minimal.

        \begin{defn}
            Set $\K=\{1\leq\xi\leq\kappa:\xi \text{ a cardinal}\}$;
            then $|\K|=\max\{|\alpha|,\aleph_0\}$. Define for every
            $f\in\Oo$ a function
            $$
                s_f:\quad\begin{matrix} \K &\To& \K\cup\{0\}\\
                \xi&\mapsto& |\{y\in X: |f\inv[y]|=\xi\}| \end{matrix}
            $$
            In words, the function assigns to every $1\leq\xi\leq\kappa$ the
            number of equivalence classes in the kernel of $f$ which have cardinality $\xi$. We call $s_f$
            the \emph{kernel sequence} of $f$. The \emph{support}
            $\supp(s_f)$ of $s_f$ is the set of all $\xi\leq\kappa$ for which
            $s_f(\xi)\neq 0$. The \emph{strong support} $\supp'(s_f)$ of $s_f$ is the set
            of those cardinals $\xi\leq\kappa$ for which
            $s_f(\xi)\cdot\xi > |X\sm f[X]|$. The \emph{weak support} of $s_f$ is
            defined to equal $\supp(s_f)\sm\supp'(s_f)$. The restriction of
            $s_f$ to its strong support is denoted by $s_f'$. We
            write $s_f'=s_g'$ iff $\supp'(s_f)=\supp'(s_g)$ and
            $s_f$ and $s_g$ agree on $\supp'(s_f)$. For $\psi_1,\psi_2\leq\kappa$ we set
            $s_f(>\psi_1)=\sum_{\psi_1 <\zeta\leq\kappa}
            s_f(\zeta)$, and $s_f(>\psi_1,<\psi_2)=\sum_{\psi_1 <\zeta<\psi_2}
            s_f(\zeta)$, and similarly with $\leq$ and $\geq$.
        \end{defn}

        \begin{defn}\label{DEF:characteristicValues}
            For $f\in\Oo$ we define the following cardinals:
            \begin{enumerate}
                \item $\mu_f=\min\supp(s_f)$.
                \item $\sigma_f=s_f(\mu_f)$.
                \item $\varrho_f=s_f(>\mu_f)$.
                \item $\nu_f=|X\sm f[X]|$.
                \item $\varepsilon_f=\sup\supp(s_f)$.
                \item $\varepsilon_f'=\sup\supp'(s_f)$, if
                $\supp'(s_f)\neq\emptyset$, and $\efp=\mu_f$
                otherwise.
                \item $\lambda_f'=\sup\{\xi\in \supp'(s_f):\xi\leq\nu_f\}$, if that set is non-void,
                and $\lambda_f'=\mu_f$ otherwise.
                \item $\chi_f=\min\{1\leq\xi\leq\kappa: \exists
                \zeta\in\supp(s_f): s_f(\geq\xi)\leq\zeta\}$.
            \end{enumerate}
        \end{defn}

        So $1\leq\mu_f\leq\lambda_f'\leq\efp\leq\ef\leq\kappa$, and $\supp(s_f)$
        is a subset of the interval $[\mu_f,\ef]$ and
        $\supp'(s_f)$ one of the interval $[\mu_f,\efp]$. The size
        of the complement $\nu_f$ is independent of the
        other cardinals, and it will be important in our proof
        whether or not $\ef>\nu_f$, that is, whether or not there exists a kernel class larger than the complement
        of the range of $f$. If $\ef>\nu_f$, then
        $\efp=\ef$ so we can forget about $\efp$ and have either
        $\mu_f\leq\lambda_f'\leq\nu_f<\ef$ or
        $\nu_f<\mu_f\leq\ef$; in the latter case we left away
        $\lambda_f'$ as it equals $\mu_f$.
        If $\ef\leq\nu_f$, then we have
        $\mu_f\leq\efp\leq\ef\leq\nu_f$; $\lambda_f'$ is irrelevant as it equals
        $\efp$. In that case,
        $\chi_f$ will play a role and in the relevant situations (e.g. if $f$ is $\S$-minimal, or if it
        satisfies conditions
        ($\sigma$) and ($\chi$) of Theorem \ref{THM:mainTheorem}) we have $\efp<\chi_f\leq\ef^+$, where $\ef^+$ is the successor
        of $\ef$.

        \begin{thm}\label{THM:mainTheorem}
            The constant functions are $\S$-minimal. If $f\in\Oo\sm\S$ is
            nonconstant, then $f$ is $\S$-minimal if and only
            if all of the following hold:
            \begin{enumerate}
                \item[($\mu$)] $\mu_f=1$ or $\mu_f$ is infinite.
                \item[($\nu$)] If $\mu_f$ is finite, then $\nu_f$ is
                infinite or zero.
                \item[($\sigma$)] $\sigma_f=\kappa$.
                \item[($\rho$)] $\rho_f<\kappa$.
                \item[(s'dec)] $s_f'$ is strictly decreasing.
                \item[(n)] $n\nin \supp'(s_f)$ for all $1<n<\aleph_0$.
                \item[($\varepsilon$)] $\varepsilon_f=1$ or $\varepsilon_f$ is infinite.
                \item[(scont)] $s_f(\geq\xi)=\min\{s_f(\geq\zeta):\zeta<\xi\}$ for all singular
                $\xi\leq\chi_f$ and\\ $s_f(\geq n)=s_f(\geq 2)$ for all finite $2\leq n\leq\chi_f$.
                \item[($\chi$)] If $\ef\leq\nu_f$, then $s_f(\geq\chi_f)$ is finite.
                \item[(\#$\varepsilon$)] If $\ef>\nu_f$, then $s_f(\varepsilon_f)$ is infinite.
                \item[($\lambda$')] If $\ef>\nu_f$, then $s_f(\xi)=0$ for all $\lambda_f'<\xi\leq\nu_f$.
            \end{enumerate}
        \end{thm}

        The following theorem describes the clones generated by
        $\S$-minimal functions. It says that the clone an $\S$-minimal
        function $f$ generates contains those non-permutations $g$
        which satisfy the conditions of Theorem \ref{THM:mainTheorem}, have
        the same characteristic values as $f$ as defined in Definition
        \ref{DEF:characteristicValues}, agree with $f$ on the strong
        support, have the same inversely-accumulated kernel sequence
        $s_g(\geq\xi)$ below $\chi_g$ as $f$, and for which
        $\varepsilon_g$ is obtained as a maximum of the support of $s_g$ iff
        it is a maximum of the support of $s_f$.

        \begin{thm}\label{THM:whenEqual}
            Let $f,g$ be $\S$-minimal. Then $\fwithS=\gwithS$ if and only
            if all of the following hold:
            \begin{enumerate}
                \item $\mu_g=\mu_f$
                \item $\nu_g=\nu_f$
                \item $s_g'=s_f'$
                \item $\chi_g=\chi_f$
                \item $s_g(\geq\xi)=s_f(\geq\xi)$ for all
                $\xi<\chi_f$
                \item $\varepsilon_g=\ef$
                \item $s_g(\varepsilon_g)=0$ iff
                $s_f(\varepsilon_f)=0$.
            \end{enumerate}
        \end{thm}

        \begin{cor}\label{COR:numberOfMinimal}
            The number of clones (monoids) minimal in $[\cl{\S},\O]$ (in $[\S,\Oo]$)  on an
            infinite set of cardinality $\aleph_\alpha$ is
            $\max\{|\alpha|,\aleph_0\}$.
        \end{cor}

        Let $X$ be countably infinite. For all $\nu<\aleph_0$, define a monoid
        $\I_\nu$ to consist of $\S$ plus all functions $f\in\Oo$ with
        $\mu_f=\aleph_0$ and
        $\nu_f=\nu$. Denote by $\H$ the monoid containing $\S$ and all
        functions with $\ef=1$ and $\nu_f=\aleph_0$, and by $\Const$ the monoid
        of all constant operations plus the permutations.

        \begin{cor}\label{COR:countableCase}
            On countably infinite $X$, the minimal monoids above
            $\S$ are exactly the monoids
            $\I_\nu$ ($\nu<\aleph_0$), $\Const$, and $\H$.
        \end{cor}

        \subsection{Notation and notions}
            The smallest clone containing a set of functions
            $\F\subseteq\O$ is denoted by $\cl{\F}$.
            For $\xi\leq\kappa$ a cardinal we set
            $Y^f_\xi=\{y\in X: |f\inv [y]|=\xi\}$ and $Y^f_{>\xi}=\bigcup_{\zeta>\xi}
            Y^f_\zeta$. Similarly we use notations like $Y^f_{\geq\xi}$ and so
            on; a less self-explanatory one is $Y^f_{>\xi,<\zeta}$ which we define to be $Y^f_{>\xi}\cap Y^f_{<\zeta}$.
            We say that a set $Y\sub X$ is \emph{large} iff
            $|Y|=\kappa$, and that it is \emph{small} otherwise.
            $Y$ is \emph{co-large} (\emph{co-small}) iff its
            complement in $X$ is large (small). For a function
            $f\in\Oo$, we denote the image of $Y\sub X$ under $f$
            by $f[Y]$ and the preimage by $f\inv[Y]$; if $Y=\{y\}$
            is a singleton, then we cut short and write $f\inv[y]$
            rather than $f\inv[\{y\}]$. Since we are interested in
            cardinals as arguments and values of kernel sequences,
            a statement like ``for all $\psi_1<\xi<\psi_2$'' or
            ``for all $\xi$ in the interval $(\psi_1,\psi_2)$'' will
            usually refer to all cardinals between $\psi_1$ and
            $\psi_2$, not all ordinals; occasionally, however, we
            will enumerate a set $Z$ of cardinality $\xi$ by something like
            $Z=\{z_\zeta:\zeta<\xi\}$, in which case $\zeta$
            refers to all ordinals below $\xi$. We shall mention
            explicitly whenever this is the case.

        \subsection{A gimmick for the quest}

        \begin{lem}\label{LEM:basicEquivalence}
            If $f,g\in\Oo$ are unary functions satisfying $s_f=s_g$
            and $\nu_f=\nu_g$, then there exist $\beta,
            \gamma\in\S$ such that $f=\beta\circ g\circ\gamma$.
        \end{lem}
        \begin{proof}
            The assumption $s_f=s_g$ implies that there is $\gamma\in\S$ such that $f(x)=f(y)$ iff
            $g\circ\gamma(x)=g\circ\gamma(y)$, for all $x,y\in X$.
            Obviously, $|f[X]|=|g[X]|=|g\circ\gamma[X]|$ as $s_f=s_g$.
            Together with the fact that $|\c f[X]|=|\c g[X]|$ this
            implies that we can find $\beta\in\S$ such that
            $f[X]=\beta\circ g\circ\gamma[X]$, and even so that
            $f=\beta\circ g\circ\gamma$.
        \end{proof}

    \section{Sufficiencies for $\S$-minimality}

    We prove that the conditions of Theorem \ref{THM:mainTheorem}
    are sufficient for a function to be $\S$-minimal.

        \subsection{Things true about everybody}

            In this section we derive properties of functions
            generated by operations which satisfy all or some conditions of Theorem
            \ref{THM:mainTheorem}.

            \subsubsection{The man who wasn't there}

                \begin{lem}\label{LEM:NU:nu_g}
                    Let $f\in\Oo$ and $g\in\cl{\{f\}\cup\S}\sm\S$. If $\nu_f$ is infinite or zero, then
                    $\nu_g=\nu_f$. If $\nu_f$ is finite, then
                    $\nu_g\geq\nu_f$ is finite as well.
                \end{lem}
                \begin{proof}
                    It is enough to show that if $f,g\in\Oo$, then $\nu_f\leq
                    \nu_{f\circ g}\leq \nu_f+\nu_g$: The assertion then
                    follows by induction over complexity of terms. Since $f[X]\supseteq f\circ g[X]$
                    we have that $|X\sm f[X]|\leq|X\sm
                    f\circ g[X]|$, so $\nu_f\leq\nu_{f\circ g}$. Also, $\nu_{f\circ g}=|X\sm
                    f\circ g[X]|\leq|(X\sm f[X])\cup f[X\sm g[X]]|\leq |X\sm f[X]|+ |X\sm
                    g[X]|=\nu_f+\nu_g$ and we are done.
                \end{proof}

                \begin{lem}\label{LEM:NU:(nu)->NuStable}
                    Let $f\in\Oo$ satisfy ($\nu$), and let
                    $g\in\cl{\{f\}\cup\S}\sm\S$. Then $\nu_g=\nu_f$.
                \end{lem}
                \begin{proof}
                    If $\nu_f$ is infinite, then we can refer to Lemma
                    \ref{LEM:NU:nu_g}, so
                    assume that $\nu_f$ is finite. If $\mu_f\leq\nu_f$, then $\nu_f>0$, and we would have
                    to have that $\mu_f$ is infinite by condition ($\nu$),
                    so this case cannot occur. Assume therefore that
                    $\mu_f>\nu_f$. Using induction over terms, it is
                    enough to show that if $h\in\Oo$ satisfies $\nu_h=\nu_f$, then $\nu_{f\circ h}=\nu_f$.
                    Indeed, since $\mu_f>|X\sm f[X]|=|X\sm h[X]|$,
                    $h$ hits every class of the kernel of $f$, so that
                    $f\circ h[X]=f[X]$ and we are done.
                \end{proof}

            \subsubsection{The dwarf-box}

                \begin{lem}\label{LEM:UBOUND:FINITE:sh}
                    Let $f,g\in\Oo$, where $g$ satisfies \cnu,
                    and set $h=f\circ g$. If $1\leq n<\aleph_0$, then
                    $s_h(n)\leq s_f(n)+s_g(>1,\leq n)+\min(\nu_g,s_f(>n,\leq\nu_g))$.
                \end{lem}
                \begin{proof}
                    If $\mu_g$ is infinite, then $s_h(n)=0$ for
                    all $1\leq n<\aleph_0$ so there is nothing to
                    show. So $\mu_g$ is finite and thus $\nu_g$ is
                    zero or infinite, by \cnu.
                    If $|h\inv [y]|=n$ for some $y\in
                    X$, then $|f\inv[y]|\geq n$ or there exists $z\in
                    f\inv[y]$ with $1<|g\inv[z]|\leq n$. The latter case
                    occurs at most $s_g(>1,\leq n)$ times. $|f\inv[y]|=n$
                    occurs exactly $s_f(n)$ times. If $|f\inv[y]|> n$
                    then there exists $z\in f\inv[y]$ not in the range of
                    $g$, which happens at most $\nu_g$ times.
                    Also, in that case we cannot have $|f\inv[y]|> \nu_g$, for
                    otherwise $|h\inv[y]|\geq |f\inv[y]\cap
                    g[X]|=|f\inv[y]|$, the latter equality holding
                    as $\nu_g$ is zero or infinite. Hence,
                    $s_h(n)\leq s_f(n)+s_g(>1,\leq n)+\min(\nu_g,s_f(>n,\leq\nu_g))$.
                \end{proof}
                \begin{lem}\label{LEM:UBOUND:FINITE:sg}
                    Let $f\in\Oo$ satisfy \cnu. Then we have for all
                    $g\in\fwithS$: For all $1<
                    n<\aleph_0$, if $s_f(>1,\leq
                    n)+\min(\nu_f,s_f(>n,\leq\nu_f))$ is zero or
                    infinite, then
                    $s_g(n)\leq s_f(>1,\leq n)+\min(\nu_f,s_f(>n,\leq\nu_f))$; otherwise $s_g(n)$ is finite.
                \end{lem}
                \begin{proof}
                    We prove this by induction over terms. We can obviously assume that $\mu_f$ is finite, hence
                    $\nu_f$ is zero or infinite by \cnu. The
                    statement is clear if $g\in\{f\}\cup\S$, so
                    assume $g=f\circ h$, with $h\in\fwithS$
                    satisfying the induction hypothesis; by Lemma \ref{LEM:NU:nu_g}, $\nu_h=\nu_f$ and in particular
                    $h$ satisfies \cnu. Therefore,
                    $s_g(n)\leq s_f(n)+s_h(>1,\leq n)+\min(\nu_f,s_f(>n,\leq\nu_f))$ by
                    Lemma \ref{LEM:UBOUND:FINITE:sh}. Now observe
                    that if $1<k\leq n$, then $s_f(>1,\leq
                    k)+\min(\nu_f,s_f(>k,\leq\nu_f))\leq
                    s_f(>1,\leq n)+\min(\nu_f,s_f(>n,\leq\nu_f))$.
                    Therefore, if $s_f(>1,\leq
                    n)+\min(\nu_f,s_f(>n,\leq\nu_f))$ is finite,
                    then so is $s_f(>1,\leq
                    k)+\min(\nu_f,s_f(>k,\leq\nu_f))$, and thus $s_h(k)$ is finite by induction
                    hypothesis, for all $1<k\leq n$. Hence,
                    $s_g(n)$ is finite. If on the other hand $s_f(>1,\leq
                    n)+\min(\nu_f,s_f(>n,\leq\nu_f))$ is infinite or zero,
                    then we have $s_h(k)\leq s_f(>1,\leq
                    n)+\min(\nu_f,s_f(>n,\leq\nu_f))$ for all
                    $1<k\leq n$ by induction hypothesis, finishing the
                    proof.
                \end{proof}

                \begin{lem}\label{LEM:STABLE:(Nonu)+(nu)Stable}
                    Let $f\in\Oo$ satisfy (n) and ($\nu$), and let
                    $g\in\fwithS\sm\S$. Then $g$ satisfies (n) as well.
                \end{lem}
                \begin{proof}
                    By ($\nu$) and Lemma \ref{LEM:NU:(nu)->NuStable} we have that $\nu_g=\nu_f$.
                    If $1<\nu_f<\aleph_0$, then
                    $\mu_f$ is infinite and so is $\mu_g$, so there is nothing to
                    show. If $\nu_f$ is zero or infinite, then since $s_f(n)\leq\nu_f$ for all
                    $1<n<\aleph_0$ by (n), we have that the same holds for $s_g$ by Lemma
                    \ref{LEM:UBOUND:FINITE:sg}. Hence, $n\nin \supp'(s_g)$ for all $1<n<\aleph_0$.
                \end{proof}

            \subsubsection{Upper bounds}

                \begin{lem}\label{LEM:UPPERBOUND:sg(>xi)}
                    Let $f\in\Oo$, and let $\xi\leq\kappa$ be infinite.
                    Then we have for all $g\in\fwithS$: If $s_f(>\xi)$ is infinite or zero, then $s_g(>\xi)\leq
                    s_f(>\xi)$; if $s_f(>\xi)$ is finite, then $s_g(>\xi)$ is finite as well.
                \end{lem}
                \begin{proof}
                    Using induction over terms, it is sufficient to show
                    that if $f,g\in\Oo$, then $s_{f\circ g}(>\xi)\leq
                    s_f(>\xi)+s_g(>\xi)$. Indeed, let $y\in X$ with $|(f\circ g)\inv [y]|>\xi$. Then $|f\inv
                    [y]|>\xi$ or there exists $z\in f\inv [y]$ with $|g\inv
                    [z]|>\xi$. The first possibility occurs for
                    $s_f(>\xi)$ elements $y\in X$, and the second one for
                    $s_g(>\xi)$ elements $y\in X$ and we are done.
                \end{proof}

                \begin{lem}\label{LEM:UBOUND:sh(xi)}
                    Let $f,g\in\Oo$, and set $h=f\circ g$. Let
                    $\xi\leq \kappa$ be infinite and regular. Then
                    $s_h(\xi)\leq s_f(\xi)+s_g(\xi)+\min(\nu_g,s_f(>\xi,\leq\nu_g))$.
                \end{lem}
                \begin{proof}
                    If $|h\inv[y]|=\xi$, then either there exists $z\in
                    f\inv[y]$ with $|g\inv[z]|=\xi$, or $|f\inv[y]|\geq\xi$, because
                    $\xi=\sum_{z\in f\inv[y]}|g\inv[z]|$ is infinite and regular. The first case can occur at most $s_g(\xi)$
                    times. That $|f\inv[y]|=\xi$ occurs $s_f(\xi)$ times,
                    so let us consider the last possibility,
                    $|f\inv[y]|>\xi$. If $|f\inv[y]|>\nu_g$, then
                    $|h\inv[y]|\geq |f\inv[y]\cap g[X]|=|f\inv [y]|>\xi$, so
                    we must have
                    $\xi <|f\inv[y]|\leq \nu_g$. Only
                    $s_f(>\xi,\leq\nu_g)$ elements $y\in X$ have this
                    property. Moreover, if $|f\inv[y]|>\xi$ but $|h\inv[y]|=\xi$,
                    then there exists $z\in f\inv[y]\sm g[X]$,
                    which happens at most $\nu_g$ times.
                \end{proof}

                \begin{lem}\label{LEM:UBOUND:regular}
                    Let $f\in\Oo$, and let $\xi\leq\kappa$ be infinite and regular. Then
                    for all $g\in\cl{\{f\}\cup\S}$ we have: If $s_f(\xi)+\min(\nu_f,s_f(>\xi,\leq\nu_f))$ is infinite or
                    zero, then $s_g(\xi)\leq s_f(\xi)+\min(\nu_f,s_f(>\xi,\leq\nu_f))$. If
                    $s_f(\xi)+\min(\nu_f,s_f(>\xi,\leq\nu_f))$ is finite,
                    then $s_g(\xi)$ is finite as well.
                \end{lem}
                \begin{proof}
                    We use induction over terms.
                    The lemma is clear if $g=f$, so say
                    $g=t\circ q$, with $q\in\fwithS$ satisfying the induction
                    hypothesis, and $t\in \{f\}\cup\S$. There is nothing
                    to show if $t\in\S$ so say $t=f$. By Lemma
                    \ref{LEM:UBOUND:sh(xi)} we have $s_g(\xi)\leq
                    s_f(\xi)+s_q(\xi)+\min(\nu_q,s_f(>\xi,\leq\nu_q))$.
                    We distinguish two
                    cases:\\
                    \textbf{Case 1}.
                    If $\nu_f$ is infinite, then $\nu_q=\nu_f$
                    by Lemma \ref{LEM:NU:nu_g}, and thus
                    $s_g(\xi)\leq
                    s_f(\xi)+s_q(\xi)+\min(\nu_f,s_f(>\xi,\leq\nu_f))$.
                    Now if $s_f(\xi)+\min(\nu_f,s_f(>\xi,\leq\nu_f))$ is
                    infinite or zero, then using the induction hypothesis for $q$ we get $s_g(\xi)\leq
                    2\cdot (s_f(\xi)+\min(\nu_f,s_f(>\xi,\leq\nu_f)))$;
                    the assertion clearly follows. If $s_f(\xi)+\min(\nu_f,s_f(>\xi,\leq\nu_f))$ is
                    finite, then $s_q(\xi)$ is finite too and so is
                    $s_g(\xi)$.\\
                    \textbf{Case 2}. If $\nu_f$ is finite, then so is
                    $\nu_q$ by Lemma \ref{LEM:NU:nu_g}, so $s_g(\xi)\leq
                    s_f(\xi)+s_q(\xi)+\min(\nu_q,s_f(>\xi,\leq\nu_q))\leq s_f(\xi)+
                    s_q(\xi)$ as $\xi>\nu_q$. Now if $s_f(\xi)$ is infinite, then
                    $s_q(\xi)\leq s_f(\xi)+\min(\nu_f,s_f(>\xi,\leq\nu_f))=s_f(\xi)$ by induction
                    hypothesis, so $s_g(\xi)\leq s_f(\xi)$. If $s_f(\xi)$
                    is finite, then so is $s_q(\xi)$ by induction
                    hypothesis; hence, $s_g(\xi)\leq s_f(\xi)+s_q(\xi)$ is
                    finite.
                \end{proof}

                \begin{lem}\label{LEM:UBOUND:sh(xi)singular}
                    Let $f,g\in\Oo$, and set $h=f\circ g$. Let
                    $\xi\leq\kappa$ be infinite, and let $\lambda < \xi$. Then
                    $s_h(\xi)\leq s_f(\xi)+ s_g(>\lambda,\leq \xi)+\min(\nu_g, s_f(>\xi,\leq\nu_g))$.
                \end{lem}
                \begin{proof}
                    If $|h\inv[y]|=\sum_{z\in f\inv[y]}|g\inv[z]|=\xi$, then either there exists $z\in
                    f\inv[y]$ with $\lambda< |g\inv[z]|\leq\xi$, or $|f\inv[y]|\geq\xi$.
                    The first case can occur at most $s_g(>\lambda,\leq \xi)$
                    times. In the second case, observe that if
                    $|f\inv[y]|>\xi$, then there must exist $z\in
                    f\inv[y]$ which is not in the range of $g$; this can
                    happen at most $\nu_g$ times. Also, in that case we
                    must have $|f\inv[y]|\leq\nu_g$, for otherwise
                    $|h\inv[y]|\geq|f\inv[y]\cap g[X]|=|f\inv[y]|>\xi$; the two conditions for $f$
                    are satisfied by $\min(\nu_g,s_f(>\xi,\leq\nu_g))$ elements
                    $y\in X$. The last possibility is that
                    $|f\inv[y]|=\xi$, which happens at most $s_f(\xi)$
                    times.
                \end{proof}

                \begin{lem}\label{LEM:UPPERBOUND:sg(>xi,<nu)}
                    Let $f\in\Oo$, and let $\xi\leq\kappa$ be infinite.
                    Then we have for all $g\in\fwithS$:
                    If $s_f(>\xi,\leq\nu_f)$ is infinite or zero, then $s_g(>\xi,\leq\nu_f)\leq
                    s_f(>\xi,\leq\nu_f)$; if $s_f(>\xi,\leq\nu_f)$ is finite,
                    then $s_g(>\xi,\leq\nu_f)$ is finite as well.
                \end{lem}
                \begin{proof}
                    Using induction over terms, it is sufficient to show
                    that if $h\in\Oo$, then $s_{h\circ f}(>\xi,\leq\nu_f)\leq
                    s_h(>\xi,\leq\nu_f)+s_f(>\xi,\leq\nu_f)$. Indeed, let
                    $y\in X$ with $\xi<|(h\circ f)\inv [y]|\leq\nu_f$. Then $|h\inv
                    [y]|>\xi$ or there exists $z\in h\inv [y]$ with $\xi<|f\inv
                    [z]|\leq\nu_f$. The latter possibility occurs for at
                    most $s_f(>\xi,\leq\nu_f)$ elements $y\in X$. If $|h\inv
                    [y]|>\xi$ and $|h\inv[y]|>\nu_f$, then $|(h\circ
                    f)\inv[y]|\geq |h\inv[y]\cap f[X]|=|h\inv[y]|>\nu_f$, so this is impossible and we have to have
                    $|h\inv[y]|\leq\nu_f$. Therefore, the first case
                    happens at most $s_h(>\xi,\leq\nu_f)$ times.
                \end{proof}

                \begin{lem}\label{LEM:UBOUND:singular}
                    Let $f\in\Oo$, and let $\xi\leq\kappa$ be infinite. Let moreover $\lambda <\xi$. Then
                    for all $g\in\cl{\{f\}\cup\S}$ we have: If $s_f(>\lambda,\leq \xi)+\min(\nu_f,
                    s_f(>\xi,\leq\nu_f))$ is infinite or zero, then
                    $s_g(\xi)\leq s_f(>\lambda,\leq \xi)+\min(\nu_f,
                    s_f(>\xi,\leq\nu_f))$. If $s_f(>\lambda,\leq \xi)+\min(\nu_f,
                    s_f(>\xi,\leq\nu_f))$ is finite, then $s_g(\xi)$ is
                    finite as well.
                \end{lem}
                \begin{proof}
                    We use induction over terms.
                    The lemma is clear if $g=f$, so say
                    $g=q\circ t$, with $q\in\fwithS$ satisfying the induction
                    hypothesis, and $t\in \{f\}\cup\S$. There is nothing
                    to show if $t\in\S$ so say $t=f$. By Lemma
                    \ref{LEM:UBOUND:sh(xi)singular}, we have $s_g(\xi)\leq
                    s_q(\xi)+ s_f(>\lambda,\leq \xi)+\min(\nu_f,
                    s_q(>\xi,\leq\nu_f))$. We distinguish three cases:\\
                    \textbf{Case 1.} Assume first $s_f(>\lambda,\leq \xi)+\min(\nu_f,
                    s_f(>\xi,\leq\nu_f))=0$; we have to show $s_g(\xi)=0$. We have $s_q(\xi)=0$ by induction hypothesis, so
                    $s_g(\xi)\leq s_f(>\lambda,\leq \xi)+\min(\nu_f,s_q(>\xi,\leq\nu_f))=\min(\nu_f,s_q(>\xi,\leq\nu_f))$.
                    Now if $\nu_f=0$, then we have $s_g(\xi)=0$ and we are
                    done. If $\nu_f>0$, then $s_f(>\xi,\leq\nu_f)=0$ since
                    $\min(\nu_f,
                    s_f(>\xi,\leq\nu_f))=0$, and Lemma
                    \ref{LEM:UPPERBOUND:sg(>xi,<nu)} implies
                    $s_q(>\xi,\leq\nu_f)=0$, so $s_g(\xi)=0$.\\
                    \textbf{Case 2.} Now assume that $s_f(>\lambda,\leq \xi)+\min(\nu_f,
                    s_f(>\xi,\leq\nu_f))$ is infinite. Then
                    $s_q(\xi)\leq s_f(>\lambda,\leq \xi)+\min(\nu_f,
                    s_f(>\xi,\leq\nu_f))$ by induction hypothesis, so
                    $s_g(\xi)\leq 2\cdot s_f(>\lambda,\leq \xi)+\min(\nu_f,
                    s_f(>\xi,\leq\nu_f))+\min(\nu_f,
                    s_q(>\xi,\leq\nu_f))$. If $s_f(>\xi,\leq\nu_f)$ is infinite, then by Lemma
                    \ref{LEM:UPPERBOUND:sg(>xi,<nu)} we have
                    $s_q(>\xi,\leq\nu_f)\leq s_f(>\xi,\leq\nu_f)$, so
                    $s_g(\xi)\leq 2\cdot (s_f(>\lambda,\leq \xi)+\min(\nu_f,
                    s_f(>\xi,\leq\nu_f)))$ and we are done. If
                    $s_f(>\xi,\leq\nu_f)$ is finite, then
                    $s_q(>\xi,\leq\nu_f)$ is finite as well by Lemma
                    \ref{LEM:UPPERBOUND:sg(>xi,<nu)}. Also in that case, $s_f(>\lambda,\leq
                    \xi)$ must be infinite, so $s_g(\xi)\leq 2\cdot s_f(>\lambda,\leq \xi)+\min(\nu_f,
                    s_f(>\xi,\leq\nu_f))+\min(\nu_f,
                    s_q(>\xi,\leq\nu_f))= 2\cdot s_f(>\lambda,\leq \xi)=s_f(>\lambda,\leq
                    \xi)$.\\
                    \textbf{Case 3.} We consider the case where $s_f(>\lambda,\leq \xi)+\min(\nu_f,
                    s_f(>\xi,\leq\nu_f))$ is finite. By induction
                    hypothesis, $s_q(\xi)$ is finite. Now if $\nu_f$ is finite,
                    then $s_g(\xi)\leq
                    s_q(\xi)+ s_f(>\lambda,\leq \xi)+\min(\nu_f,
                    s_q(>\xi,\leq\nu_f))$ is finite, too. If $\nu_f$ is
                    infinite, then $s_f(>\xi,\leq\nu_f)$ must be finite,
                    and thus $s_q(>\xi,\leq\nu_f)$ is finite by Lemma
                    \ref{LEM:UPPERBOUND:sg(>xi,<nu)}; again, we have that
                    $s_g(\xi)\leq
                    s_q(\xi)+ s_f(>\lambda,\leq \xi)+\min(\nu_f,
                    s_q(>\xi,\leq\nu_f))$ is finite.
                \end{proof}

            \subsubsection{Lower bounds}

                \begin{lem}\label{LEM:LBOUND:BEYONDNU:sh(xi)}
                    Let $f,g\in\Oo$, and set $h=f\circ g$. Let
                    $\xi\in (\nu_g,\kappa]$ be infinite, and assume
                     that either $s_g(>\xi)=0$ or $s_g(>\xi)<s_f(\xi)$ and $s_f(\xi)$ is infinite.
                    Then $s_h(\xi)\geq s_f(\xi)$.
                \end{lem}
                \begin{proof}
                    Fix some
                    $y\in Y^f_\xi$. Then $|h\inv [y]|\geq |f\inv [y]\cap
                    g[X]|$. Since $\nu_g <\xi$ we have that $|f\inv [y]\cap
                    g[X]|=|f\inv[y]|=\xi$. Now $|h\inv [y]| >\xi$ if and only if
                    there exists $z\in f\inv[y]$ with $|g\inv[z]|>\xi$.
                    This happens only for $s_g(>\xi)$ elements
                    $y\in Y^f_\xi$, so it does not happen for $s_f(\xi)$ elements $y\in Y^f_\xi$, since
                    either $s_g(>\xi)=0$ or $s_f(\xi)$ is infinite and $s_g(>\xi)<s_f(\xi)$. Hence,
                    $|h\inv[y]|=\xi$ for at least $s_f(\xi)$ elements
                    $y\in Y^f_\xi$.
                \end{proof}

                \begin{lem}\label{LEM:LBOUND:BEYONDNU:sg(xi)lowerbound}
                    Let $f\in\Oo$, let $\xi\in(\nu_f,\kappa]$ be infinite, and
                    assume that either $s_f(>\xi)=0$ or
                    $s_f(>\xi)<s_f(\xi)$ and $s_f(\xi)$ is infinite. Then
                    for all $g\in\cl{\{f\}\cup\S}\sm\S$ we have $s_g(\xi)\geq
                    s_f(\xi)$.
                \end{lem}
                \begin{proof}
                    Because $\xi>\nu_f$ is infinite, we have $\xi>\nu_g$ for all $g\in\fwithS$ by Lemma
                    \ref{LEM:NU:nu_g}. If $s_f(>\xi)=0$, then $s_g(>\xi)=0$ for all $g\in\fwithS$ by Lemma
                    \ref{LEM:UPPERBOUND:sg(>xi)}. Also, by the same lemma, we
                    have that if $s_f(\xi)$ is infinite, then
                    $s_f(>\xi)<s_f(\xi)$ implies $s_g(>\xi)<s_f(\xi)$ for all
                    $g\in\fwithS$. The rest of the proof is induction over terms and
                    Lemma \ref{LEM:LBOUND:BEYONDNU:sh(xi)}.
                \end{proof}

                \begin{lem}\label{LEM:LBOUND:BELOWNU:sh(xi)}
                    Let $f,g\in\Oo$, and set $h=f\circ g$. Let
                    $\xi\leq\nu_g$ be infinite, and assume $s_g(>\xi)+\nu_g <s_f(\xi)$
                    and that $s_f(\xi)$ is infinite.
                    Then $s_h(\xi)\geq s_f(\xi)$.
                \end{lem}
                \begin{proof}
                    Fix some $y\in Y^f_\xi$. Then $|h\inv [y]|\geq |f\inv [y]\cap
                    g[X]|$. We have that $|f\inv [y]\cap
                    g[X]|\geq\xi$ for at least $s_f(\xi)$ elements $y\in Y^f_\xi$, since $|f\inv [y]\cap
                    g[X]|<\xi$ implies that $f\inv[y]\sm g[X]$ is non-empty, which happens at most $\nu_g<s_f(\xi)$ times.
                    Now $|h\inv [y]| >\xi$ if and only if
                    there exists $z\in f\inv[y]$ with $|g\inv[z]|>\xi$.
                    By assumption $s_g(>\xi)<s_f(\xi)$ this happens for fewer than $s_f(\xi)$ elements
                    $y\in Y^f_\xi$, so that
                    $|h\inv[y]|=\xi$ for $s_f(\xi)$ elements
                    $y\in Y^f_\xi$.
                \end{proof}

                \begin{lem}\label{LEM:LBOUND:BELOWNU:lowerbound}
                    Let $f\in\Oo$, let $\xi\leq\nu_f$ be infinite,
                    and assume $s_f(>\xi)+\nu_f <s_f(\xi)$ and that $s_f(\xi)$ is
                    infinite. Then
                    for all $g\in\cl{\{f\}\cup\S}\sm\S$ we have $s_g(\xi)\geq
                    s_f(\xi)$.
                \end{lem}
                \begin{proof}
                    Because $s_f(\xi)>\nu_f$ is infinite, we have $s_f(\xi)>\nu_g$ for all $g\in\fwithS$ by Lemma
                    \ref{LEM:NU:nu_g}. By Lemma \ref{LEM:UPPERBOUND:sg(>xi)}, since
                    $s_f(>\xi)<s_f(\xi)$,
                    and since $s_f(\xi)$ is infinite, we
                    have $s_g(>\xi)<s_f(\xi)$ for all
                    $g\in\fwithS$. The rest of the proof is induction over terms and
                    Lemma \ref{LEM:LBOUND:BELOWNU:sh(xi)}.
                \end{proof}

            \subsubsection{The king}

                \begin{lem}\label{LEM:EPSILON:eps=1OrInfinite->Stable}
                    Let $f\in\Oo$ satisfy $(\varepsilon$).
                    Then $\varepsilon_g=\varepsilon_f$ for all
                    $g\in\fwithS\sm\S$.
                \end{lem}
                \begin{proof}
                    If $\ef=1$, then
                    $f$ is injective and so is $g$, hence
                    $\varepsilon_g=1$. Otherwise $\ef$ is infinite.
                    Fix $\xi\leq\kappa$; if $s_f(\geq\xi)>0$, then clearly also
                    $s_g(\geq\xi)>0$ since kernel classes cannot become
                    smaller, so $\varepsilon_g\geq \ef$.  On the other hand,
                    $s_g(>\varepsilon_f)\leq s_f(>\ef)=0$ by Lemma \ref{LEM:UPPERBOUND:sg(>xi)}, so
                    $\varepsilon_g\leq\ef$.
                \end{proof}

                \begin{defn}
                    We say that $f\in\Oo$ satisfies ($\varepsilon$reg) iff
                    $s_f(\varepsilon_f)>0$ or $\ef$ is regular.
                \end{defn}
                \begin{lem}\label{LEM:SUF:(epsreg)}
                    If $f\in\Oo$ satisfies (s'dec), (scont), and \cchi, then it satisfies ($\varepsilon$reg).
                \end{lem}
                \begin{proof}
                    If $\ef>\nu_f$, then the support of $s_f$ above $\nu_f$ is finite by
                    (s'dec), so $s_f(\ef)>0$. If $\ef\leq\nu_f$ and
                    $\ef>\chi_f$, then $0<s_f(\geq\chi_f)<\aleph_0$ by
                    \cchi, so again $s_f(\ef)>0$. If $\ef\leq\nu_f$ and
                    $\ef\leq\chi_f$, then we have that if $\ef$ is
                    singular, then
                    $s_f(\ef)=s_f(\geq\ef)=\min\{s_f(\geq\zeta):\zeta<\ef\}>0$,
                    by (scont).
                \end{proof}

            \subsubsection{Farmers}

                \begin{defn}
                    We say that $f\in\Oo$ satisfies ($\kappa$) iff
                    $\nu_f=\kappa$ implies $s_f(\kappa)=0$.
                \end{defn}
                \begin{lem}\label{LEM:SUF:(kappa)}
                    If $f\in\Oo$ satisfies \csig and \cchi, then it satisfies ($\kappa$).
                \end{lem}
                \begin{proof}
                    Observe that $\nu_f=\kappa$ implies $\ef\leq\nu_f$. If $s_f(\kappa)>0$, then $\chi_f=\mu_f$. But then
                    $s_f(\geq\chi_f)=s_f(\geq \mu_f)\geq\sigma_f=\kappa$ by \csig,
                    contradicting \cchi.
                \end{proof}

                \begin{lem}\label{LEM:STABLE:MU}
                    Let $f\in\Oo$ satisfy ($\mu$), ($\nu$), ($\sigma$), ($\rho$), ($\varepsilon$), ($\varepsilon$reg), and
                    ($\kappa$). Then
                    $g$ satisfies $\mu_g=\mu_f$, $\nu_g=\nu_f$, $\varepsilon_g=\ef$, and ($\mu$), ($\nu$),
                    ($\sigma$), ($\rho$), ($\varepsilon$), ($\varepsilon$reg), and ($\kappa$), for all $g\in\fwithS\sm\S$.
                \end{lem}
                \begin{proof}
                    Using induction
                    over terms, we assume $g=f\circ h$, with
                    $h\in\fwithS$ having all asserted properties. We are going to
                    prove $\mu_g=\mu_f$, $\nu_g=\nu_f$,
                    $\varepsilon_g=\varepsilon_f$, ($\sigma$), ($\rho$), ($\varepsilon$), ($\varepsilon$reg), and
                    ($\kappa$); conditions ($\mu$) and ($\nu$) will follow
                    automatically from $\mu_g=\mu_f$ and $\nu_g=\nu_f$. By Lemmas \ref{LEM:NU:(nu)->NuStable} and
                    \ref{LEM:EPSILON:eps=1OrInfinite->Stable} and conditions
                    ($\nu$) and ($\varepsilon$) we have
                    $\nu_g=\nu_f$ and $\varepsilon_g=\varepsilon_f$. We prove ($\varepsilon$reg). If
                    $\varepsilon_g=\varepsilon_h$ is singular, then there exists
                    $y\in Y^h_{\varepsilon_h}$, by ($\varepsilon$reg);
                    but then $f(y)\in Y^g_{\varepsilon_g}$, and hence $g$
                    satisfies ($\varepsilon$reg). We show ($\kappa$). If $\nu_g=\nu_f=\kappa$, then
                    $s_f(\kappa)=0$ by ($\kappa$).
                    Now if $\ef<\kappa$, then
                    $\varepsilon_g=\ef<\kappa$, so $s_g(\kappa)=0$. If
                    $\ef=\kappa$, then $\kappa$ is regular by
                    ($\varepsilon$reg). But then
                    $s_g(\kappa)\leq s_f(\kappa)=0$, by
                    Lemma \ref{LEM:UBOUND:regular}. Hence, $g$ satisfies
                    ($\kappa$).\\
                    We now claim that $g[X]$ is large. Indeed, this is trivial if
                    $\nu_g<\kappa$, so assume $\nu_g=\nu_f=\kappa$. Then by ($\kappa$) we have
                    $s_g(\kappa)=s_f(\kappa)=0$. Since
                    $|\bigcup_{y\in g[X]} g\inv[y]|=\kappa$,
                    $\varepsilon_g<\kappa$ immediately
                    implies $|g[X]|=\kappa$, so consider the case
                    $\varepsilon_g=\ef=\kappa$. Because $s_f(\kappa)=0$
                    and by ($\varepsilon$reg) we have that $\kappa$ is regular.
                    Hence, $|\bigcup_{y\in g[X]} g\inv[y]|=\kappa$ again
                    implies $|g[X]|=\kappa$.\\
                    Now if $g\inv[y]>\mu_f$, then by ($\mu$)
                    we have that
                    $y\in Y^f_{>\mu_f}$ or there exists $z\in f\inv[y]\cap
                    Y^h_{>\mu_f}$, which happens for at most
                    $\rho_f+\rho_h<\kappa$ times. Therefore $Y^g_{>\mu_f}$
                    is small, and hence $Y^g_{\mu_f}$ must be large since
                    $g[X]=Y^g_{\mu_f}\cup Y^g_{>\mu_f}$ is large. Whence, $\mu_g=\mu_f$ and $g$
                    satisfies ($\sigma$) and ($\rho$).
                \end{proof}

            \subsubsection{The valley of giants}

                \begin{defn}
                    We say that $f\in\Oo$ satisfies (s'inf) iff
                    $s_f'(\xi)$ is infinite for all $\xi\in \supp'(s_f)$.
                \end{defn}

                \begin{lem}\label{LEM:SUF:(s'inf)}
                    If $f\in\Oo$ satisfies ($\nu$), (s'dec), and
                    (\#$\varepsilon$),
                    then it satisfies (s'inf).
                \end{lem}
                \begin{proof}
                    Let $\xi\in\supp'(s_f)$. If $\xi>\nu_f$, then
                    $s_f(\xi)\geq s_f(\varepsilon_f)\geq\aleph_0$, by (s'dec)
                    and (\#$\varepsilon$). If $\xi\leq\nu_f$, then
                    $s_f(\xi)>\nu_f$. If $\nu_f$ was finite, then we would have
                    $1\leq \xi\leq\nu_f<\aleph_0$, so in particular $\mu_f$ would be finite and $0<\nu_f<\aleph_0$,
                    contradicting ($\nu$).
                    Hence $\nu_f$ and thus also $s_f(\xi)$ are infinite.
                \end{proof}

                \begin{lem}\label{LEM:STRSUPP:s'dec->stable}
                    Let $f\in\Oo$ satisfy ($\mu$), ($\nu$),
                    ($\sigma$), ($\rho$), (s'dec), (n), ($\varepsilon$), (s'inf), ($\varepsilon$reg), and
                    ($\kappa$).
                    Then $s_g'=s_f'$ for all $g\in\fwithS\sm\S$.
                \end{lem}
                \begin{proof}
                    By Lemma \ref{LEM:STABLE:MU}, we have that
                    $\mu_g=\mu_f$, $\sigma_g=\sigma_f=\kappa$, $\varepsilon_g=\ef$, and
                    $\nu_g=\nu_f$. Let $\xi\in\supp'(s_f)$ so that $\xi>1$; then $\xi$ is infinite by
                    (n). Choose $\lambda<\xi$ such that
                    $\supp'(s_f)\cap(\lambda,\xi)$ is empty; this is
                    possible since the strong support is finite by condition
                    (s'dec). Also, if $\xi>\nu_f$, then we can choose
                    $\lambda>\nu_f$.
                    By Lemma \ref{LEM:UBOUND:singular} and since $s_f(\xi)$ is infinite by (s'inf)
                    we have $s_g(\xi)\leq s_f(>\lambda,\leq \xi)+\min(\nu_f,
                    s_f(>\xi,\leq\nu_f))$. But the latter expression equals
                    $s_f(\xi)$,
                    since if $\xi>\nu_f$, then $s_f(>\lambda,<\xi)=0$ by
                    the choice of $\lambda$ and also $s_f(>\xi,\leq\nu_f)=0$,
                    and if $\xi\leq\nu_f$, then
                    $s_f(>\lambda,<\xi)\leq\nu_f$ so that the equality
                    follows from the fact that $s_f(\xi)>\nu_f$. Hence,
                    $s_g(\xi)\leq s_f(\xi)$.
                    By (s'dec), (s'inf) and
                    Lemmas \ref{LEM:LBOUND:BELOWNU:lowerbound} and
                    \ref{LEM:LBOUND:BEYONDNU:sg(xi)lowerbound} we have
                    $s_g(\xi)\geq s_f(\xi)$, so $s_g(\xi)=s_f(\xi)$. If
                    $1\in\supp'(s_f)$, then $\mu_f=1$ and so
                    $s_g(1)=s_f(1)=\kappa$, since $\mu_g=\mu_f=1$ and
                    $\sigma_g=\sigma_f=\kappa$.\\
                    Now
                    let $\xi\nin \supp'(s_f)$; then $\xi\leq\nu_f=\nu_g$. Consider first the case where $\xi$ is infinite,
                    and choose $\lambda$ as before.  If $s_f(>\lambda,\leq\xi)+\min(\nu_f,
                    s_f(>\xi,\leq\nu_f))$ is infinite, then Lemma \ref{LEM:UBOUND:singular} implies
                    $s_g(\xi)\leq s_f(>\lambda,\leq \xi)+\min(\nu_f,
                    s_f(>\xi,\leq\nu_f))\leq\nu_f$, and thus $\xi
                    \nin\supp'(s_g)$. If on the other hand $s_f(>\lambda,\leq\xi)+\min(\nu_f,
                    s_f(>\xi,\leq\nu_f))$ is finite, then the same holds
                    for $s_g(\xi)$ and hence $s_g(\xi)\leq\nu_f$ as
                    $\nu_f\geq\xi$ is infinite; again, $\xi
                    \nin\supp'(s_g)$. Now consider the case where $\xi>1$ is finite. Then
                    conditions
                    (n) and ($\nu$) together with Lemma \ref{LEM:STABLE:(Nonu)+(nu)Stable} guarantee
                    that $\xi\nin\supp'(s_g)$. Finally, assume $\xi=1$. If $\xi$ is in the weak support of $s_f$, then
                    $\mu_f=\xi=1$ and we have $s_g(\xi)=s_f(\xi)=\kappa$.
                    Because $1\nin\supp'(s_f)$ we must have
                    $\nu_f=\nu_g=\kappa$, so $1\nin\supp'(s_g)$.
                    If $\xi=1\nin\supp(s_f)$, then $1\nin\supp(s_g)$, so in particular $1\nin\supp'(s_g)$.\\
                    Therefore, we have shown that $\supp'(s_g)=\supp'(s_f)$, and that
                    $s_g(\xi)=s_f(\xi)$ for all $\xi\in\supp'(s_f)$.
                \end{proof}

        \subsection{Gambling back the loss}

            We investigate which functions are are generated by operations
            satisfying (some of) the conditions of Theorem
            \ref{THM:mainTheorem}, the ultimate goal being to show
            that if $g\in\fwithS\sm\S$, where $f\in\Oo$ satisfies
            the conditions of Theorem \ref{THM:mainTheorem}, then
            $f\in\gwithS$, proving $\S$-minimality.

            \subsubsection{When the king is larger than the man who wasn't there}

                We modify functions $f\in\Oo$ below $\efp$.
                This finishes the case $\ef>\nu_f$, since in that case $\efp=\ef$.

                \begin{lem}\label{LEM:WSUPP:removalOfBoundedElementsLambda}
                    Let $f\in\Oo$ satisfy ($\mu$), ($\nu$), ($\sigma$), (s'dec), and (n).
                    Then there exists $g\in\fwithS$ such that $s_g(\xi)=0$ for all $\xi<\lambda_f'$ with $\xi\nin\supp'(s_f)$,
                    and $s_g(\xi)=s_f(\xi)$ for all other $\xi\leq \kappa$. In particular, there are no
                    elements below $\lambda_f'$ in the weak support of
                    $g$.
                \end{lem}
                \begin{proof}
                    We may assume that $\lambda'_f>\mu_f$, for the lemma is trivial
                    otherwise; thus, $\lambda_f'>1$ and so $\lambda_f'$ is infinite by  condition
                    (n). Also, from $\lambda_f'>\mu_f$ it follows that $\mu_f\leq\nu_f$, which together with \cnu
                    implies that $\nu_f$ is infinite. Set $Y=\bigcup \{Y^f_\zeta: \zeta<\lambda_f'\wedge
                    \zeta\nin\supp'(s_f)\}$; then $|Y|\leq\nu_f$. Observe that $Y$ does not contain $Y^f_{\mu_f}$, as
                    $\lambda_f'>\mu_f$ implies that the strong support of $s_f$ is non-empty and thus $\nu_f<\kappa$,
                    so $\mu_f\in\supp'(s_f)$ by \csig. Let $\alpha$ map all $y\in
                    Y^f_{\lambda_f'}$ into $f\inv [y]$.
                    Because $\lambda_f'$ is infinite
                    we have $|f\inv[Y^f_{\lambda_f'}]\sm\alpha [Y^f_{\lambda_f'}]|
                    =|f\inv[Y^f_{\lambda_f'}]|>\nu_f$, so we can extend
                    $\alpha$ by mapping $Y$ into $f\inv[Y^f_{\lambda_f'}]$ in such a way that $\alpha$ stays
                    injective. Now extend $\alpha$ again, mapping
                    $y$ into $f\inv[y]$ for all $y\in Y^f_{>\mu_f}$ for which $\alpha$ has not yet been
                    defined, and let $\alpha$ map a suitable part of $X\sm f[X]$ bijectively
                    onto $f\inv[Y]$. By ($\sigma$) we can choose $S\sub Y^f_{\mu_f}$ large such
                    that $Y^f_{\mu_f}\sm S$ is still large, and let
                    $\alpha$ map $S$ bijectively onto $f\inv[S]$.
                    Extend $\alpha$ to a bijection; this is possible as the domain of $\alpha$
                    is disjoint from $Y^f_{\mu_f}\sm S$ and its range is disjoint
                    from $f\inv[Y^f_{\mu_f}\sm S]$.\\
                    We calculate $|g\inv[y]|$ for all $y\in X$. If $y\in Y^f_{\mu_f}$,
                    then $|(\alpha\circ f)\inv [z]|\in\{0,\mu_f\}$ for all
                    $z\in f\inv[y]$, so $|g\inv[y]|\in\{0,\mu_f\}$ by \cmu; if
                    $y\in S\sub Y^f_{\mu_f}$, then $|g\inv[y]|=\mu_f$ as
                    $\alpha(y)\in f\inv [y]$. If $y\in Y$, then $g\inv [y]$ is empty as $f\inv[y]\sub
                    \alpha[X\sm f[X]]$. If $y\in Y^f_{>\mu_f}\sm Y$, then
                    $g\inv[y]\supseteq f\inv [y]$; also, $|(\alpha\circ
                    f)\inv[z]|\leq |f\inv[y]|$ for all $z\in f\inv[y]$,
                    so $|g\inv[y]|=|f\inv [y]|$ since $y\nin Y$ implies that $|f\inv [y]|$ is infinite.
                    Therefore we have $s_g(\xi)=0$ for
                    all $\xi<\lambda_f'$ outside the strong support of $s_f$, and
                    $s_g(\xi)=s_f(\xi)$ for all other $\xi\leq\kappa$. In
                    particular, since $\nu_g=\nu_f$ by ($\nu$) and Lemma
                    \ref{LEM:NU:(nu)->NuStable}, there are no elements
                    below $\lambda_f'$ in the weak support of $s_g$.
                \end{proof}

                \begin{lem}\label{LEM:WSUPP:BELOWLAMBDA:addingOfClasses}
                    Let $f\in\Oo$ satisfy  ($\mu$), ($\nu$), ($\sigma$), (s'dec), (n), and (s'inf).
                    Let $p\in\Oo$ be so that $s_p'=s_f'$, $\mu_p=\mu_f$ and $\nu_p=\nu_f$.
                    Then there exists
                    $g\in\fwithS$ such that
                    $s_g\rest_{[1,\lambda_f')}=s_p\rest_{[1,\lambda_f')}$, and
                    $s_g\rest_{[\lambda_f',\kappa]}=s_f\rest_{[\lambda_f',\kappa]}$.
                \end{lem}
                \begin{proof}
                    We assume that the strong support of $s_f$ below $\nu_f$ is non-void, for otherwise
                    $\lambda_f'=\mu_f'$ by definition and the the lemma is
                    trivial. For the same reason, we may assume that $\mu_f<\lambda_f'$; then $\lambda_f'$ and hence also
                    $\nu_f$ are infinite.
                    By Lemma
                    \ref{LEM:WSUPP:removalOfBoundedElementsLambda} there
                    exists $h\in\fwithS$ with the property that $s_h(\xi)=0$ for all
                    $\xi<\lambda_f'$ which are not in the strong support
                    of $s_f$ and such that $s_h(\xi)=s_f(\xi)$ for all other $\xi\leq
                    \kappa$; since $f$ satisfies ($\nu$), Lemma \ref{LEM:NU:(nu)->NuStable} implies that
                    $\nu_h=\nu_f$. Clearly, $s_h'=s_f'$ and
                    $\varepsilon_h=\varepsilon_f$. Now fix for
                    every $\xi<\lambda_f'=\lambda_h'$ outside the strong support of $s_p$ a set
                    $Z_\xi\leq Y_{\lambda_h'}^h$ such that $|Z_\xi|=s_p(\xi)$ and
                    such that $Z_{\xi_1}\cap Z_{\xi_2}=\emptyset$ whenever
                    $\xi_1\neq\xi_2$. This is possible since the sum over all
                    $s_p(\xi)$, where $\xi<\lambda_h'$ is not an element
                    of the strong support
                    of $s_p$, is at most $\nu_h<s_h(\lambda_h')$. Fix for
                    every $y\in Z_\xi$ a set $B_y\sub h\inv[y]$ with
                    $|B_y|=\xi$, and set $C_y=h\inv[y]\sm B_y$. Set $Z=\bigcup
                    \{Z_\xi:\xi\nin\supp'(s_h)\wedge \xi<\lambda_h' \}$, $B=\bigcup\{B_y:y\in
                    Z\}$, and
                    $C=\bigcup\{C_y:y\in Z\}$. Choose $S\subseteq
                    Y^h_{\mu_h}$ large such that $Y^h_{\mu_h}\sm S$ is still large, and $C'\sub X\sm h[X]$
                    with $|C'|=|C|$. Now let $\alpha$ map $S$ bijectively
                    onto $f\inv[S]\cup B$, $C'$ bijectively onto $C$, and
                    $Y^h_{\lambda_h'}$ onto $h\inv[Y^h_{\lambda_h'}]\sm (B\cup C)$; the latter can be done as
                    $|h\inv[Y^h_{\lambda_h'}]\sm (B\cup C)|=|h\inv[Y^h_{\lambda_h'}\sm Z]|=
                    |h\inv[Y^h_{\lambda_h'}]|=\lambda_h'\cdot
                    |Y^h_{\lambda_h'}|=|Y^h_{\lambda_h'}|$.
                    For all $y\in Y^h_{>\mu_h}\sm Y^h_{\lambda_h'}$ we let $\alpha(y)\in
                    f\inv[y]$. We extend $\alpha$ to a bijection and set
                    $g=h\circ\alpha\circ h$.\\
                    Now if $y\in Y^h_{\mu_h}$, then $|g\inv
                    [y]|\in\{0,\mu_h\}$, and if $y\in S$, then $|g\inv
                    [y]|=\mu_h$. If $y\in Y^h_{>\mu_h}\sm Z$, then $|g\inv
                    [y]|=|f\inv[y]|$. Indeed, (n) and the fact that $s_h$ vanishes
                    outside its strong support below $\lambda_h'$ imply that $|f\inv[y]|$ is infinite; the equation then follows
                    since $\alpha(y)\in f\inv[y]$ and $|(\alpha\circ f)\inv
                    [z]|\leq |f\inv[y]|$ for all $z\in f\inv[y]$ by
                    construction of $\alpha$.
                    If $y\in Z_\xi$ for some
                    $\mu_h<\xi<\lambda_h'$, then $|g\inv[y]|=|(\alpha\circ h)\inv[B_y\cup
                    C_y]|=|(\alpha\circ h)\inv[B_y]|=|B_y|\cdot
                    \mu_h=|B_y|=\xi$. Therefore,
                    $s_g(\mu_h)=s_h(\mu_h)=\kappa$, $s_g(\xi)=s_p(\xi)$
                    for all $\mu_h<\xi<\lambda_h'$ outside the strong
                    support of $s_f$, and $s_g(\xi)=s_h(\xi)=s_f(\xi)$ for all $\xi\geq\lambda_h'$ and
                    all $\xi\in\supp'(s_f)$.
                \end{proof}
                \begin{prop}\label{PROP:criterionWithClassBeyondNu}
                    Let $f\in\Oo$ be so that $\nu_f<\ef$. If $f$ moreover satisfies
                    ($\mu$), ($\nu$), ($\sigma$), ($\rho$), (s'dec), (n),
                    ($\varepsilon$), (scont),
                    ($\chi$), (\#$\varepsilon$), and ($\lambda$'), then it
                    is $\S$-minimal.
                \end{prop}
                \begin{rem}
                    Under those conditions, $f$ automatically satisfies ($\chi$) and
                    (scont):
                    Condition ($\chi$) is trivial as
                    $\ef\nleq\nu_f$. For (scont), observe that
                    $s_f(\geq\xi)=s_f(\geq \psi)$, where
                    $\psi=\min\{\zeta\in\supp'(s_f):\zeta>\xi\}$, if the latter set is not empty, which is the case
                    for all $\xi\leq\ef$ as
                    $\nu_f<\ef$.
                    Therefore, the function $s_f(\geq\xi)$ drops only at successor cardinals of elements of the strong
                    support, and hence only at infinite regular cardinals or at $2$, in accordance with (scont).
                \end{rem}
                \begin{proof}
                    Let $f$ satisfy all the conditions, and let
                    $g\in\fwithS\sm\S$. By Lemmas \ref{LEM:SUF:(epsreg)}, \ref{LEM:SUF:(kappa)} and \ref{LEM:SUF:(s'inf)},
                    $f$ satisfies ($\varepsilon$reg), ($\kappa$), and (s'inf). Therefore
                    we have $\mu_g=\mu_f$, $\nu_g=\nu_f$ and $\varepsilon_g=\ef$ by
                    Lemma \ref{LEM:STABLE:MU}. Moreover by the same
                    lemma, $g$ satisfies ($\mu$), ($\nu$) and
                    ($\sigma$). By Lemma \ref{LEM:STRSUPP:s'dec->stable}, $s_g'=s_f'$, in particular $g$
                    satisfies (s'dec), (s'inf),
                    and $\lambda_g'=\lambda_f'$. From Lemma \ref{LEM:STABLE:(Nonu)+(nu)Stable} and the fact that
                    $f$ satisfies (n) and ($\nu$) we infer that (n) holds for
                    $g$ as well. Therefore by Lemma
                    \ref{LEM:WSUPP:BELOWLAMBDA:addingOfClasses},
                    there exists $h\in\gwithS$ such that
                    $s_h\rest_{[1,\lambda_f')}=s_f\rest_{[1,\lambda_f')}$ and
                    $s_h\rest_{[\lambda_f',\kappa]}=s_g\rest_{[\lambda_f',\kappa]}$; thus,
                    $s_h\rest_{(\nu_f,\kappa]}=s_f\rest_{(\nu_f,\kappa]}$ as $s_g'=s_f'$.
                    Now $\supp(s_h)\cap(\lambda_f',\nu_f]$ is empty.
                    Indeed, $s_f(\xi)=0$ for all
                    $\xi\in(\lambda_f',\nu_f]$ by \clamp;
                    therefore, $s_h(\xi)=0$ for all infinite
                    $\xi\in(\lambda_f',\nu_f]$ by Lemma
                    \ref{LEM:UBOUND:singular}. If $\xi
                    \in(\lambda_f',\nu_f]$ is finite, then we must
                    have $\lambda_f'=\mu_f=1$, by (n) and \cmu.
                    Thus, $s_f$ yields constantly zero on $(1,\nu_f]$,
                    and so $s_h(\xi)=0$ by Lemma
                    \ref{LEM:UBOUND:FINITE:sg}. Hence, both $s_h$
                    and $s_f$ vanish on $(\lambda_f',\nu_f]$. Therefore,
                    $s_h=s_f$ so that since also $\nu_h=\nu_g=\nu_f$ by Lemma \ref{LEM:NU:(nu)->NuStable}, we
                    conclude $f\in\gwithS$.
                \end{proof}


            \subsubsection{Beyond the giants}

            First we show that if $f\in\Oo$ satisfies some of the conditions of Theorem \ref{THM:mainTheorem}, then
            $\chi_g=\chi_f$ and $s_g(\geq\xi)=s_f(\geq\xi)$ for all $\xi<\chi_f$ and all
            $g\in\fwithS$. Assuming $\ef\leq\nu_f$, we then modify functions $f\in\Oo$ above $\efp$ and below
            $\chi_f$.

                \begin{lem}\label{LEM:UBOUND:sg(>=xi)}
                    Let $f\in\Oo$ and $\xi\leq\kappa$ be infinite and regular or $\xi\leq 2$, and let $g\in\fwithS$.
                    If $s_f(\geq\xi)$ is infinite or zero, then $s_g(\geq\xi)\leq
                    s_f(\geq\xi)$. If
                    $s_f(\geq\xi)$ is finite, then $s_g(\geq\xi)$ is
                    finite as well.
                \end{lem}
                \begin{proof}
                    It is enough to show that if $h\in\Oo$, then $s_{f\circ h}(\geq\xi)\leq
                    s_f(\geq\xi)+s_h(\geq\xi)$; the lemma then clearly follows by induction over terms.
                    Indeed, if $|(f\circ h)\inv[y]|\geq\xi$,
                    then $|f\inv [y]|\geq\xi$ or there exists $z\in
                    f\inv[y]$ such that $|h\inv [z]|\geq\xi$, since $\xi$
                    is infinite and regular, or $\xi\leq 2$. The first possibility happens
                    $s_f(\geq\xi)$ and the second possibility
                    $s_h(\geq\xi)$ times.
                \end{proof}

                \begin{lem}\label{LEM:MIDSUPP:LBOUND:sf(>=xi)}
                    Let $f\in\Oo$ satisfy ($\varepsilon$) and ($\varepsilon$reg), and let $1\leq\xi<\chi_f$.
                    Then $s_g(\geq\xi)\geq s_f(\geq\xi)$
                    for all $g\in\fwithS\sm\S$.
                \end{lem}
                \begin{proof}
                    We can assume $\ef>1$, so $\ef$ is infinite.
                    Using induction over terms, it is enough to show
                    that if $h\in\fwithS$ satisfies $s_h(\geq\xi)\geq
                    s_f(\geq\xi)$, then also
                    $s_{h\circ f}(\geq\xi)\geq s_f(\geq\xi)$. By Lemma \ref{LEM:EPSILON:eps=1OrInfinite->Stable}, we
                    have $\varepsilon_h=\ef$. Consider
                    $Y^f_{\geq\xi}$; we claim that $|h[Y^f_{\geq\xi}]|=|Y^f_{\geq\xi}|=s_f(\geq\xi)$.
                    To see this, observe first that
                    $Y^f_{\geq\xi}=\bigcup_{y\in h[Y^f_{\geq\xi}]} (h\inv[y]\cap
                    Y^f_{\geq\xi})$. Now
                    $|h\inv[y]|<s_f(\geq\xi)$ for all $y\in X$. Indeed, otherwise we would
                    have $\ef=\varepsilon_h\geq|h\inv[y]|\geq s_f(\geq\xi)\geq\varepsilon_f$, the last equality holding
                    since $\xi<\chi_f$; thus, $s_h(\ef)>0$. But then $s_f(\ef)>0$ by
                    ($\varepsilon$reg) and Lemma \ref{LEM:UBOUND:regular}, so $\ef\in\supp(s_f)$ and $\ef=s_f(\geq\xi)$, contradicting
                    $\xi<\chi_f$. Therefore,
                    $|h\inv[y]\cap
                    Y^f_{\geq\xi}|\leq |h\inv[y]|<s_f(\geq\xi)$ for all $y\in
                    h[Y^f_{\geq\xi}]$. Thus if we had
                    $|h[Y^f_{\geq\xi}]|<s_f(\geq\xi)$, we could conclude
                    that $s_f(\geq\xi)$ is singular and the supremum of a set of cardinals of kernel classes of $h$, the
                    latter fact implying $s_f(\geq\xi)\leq\varepsilon_h=\ef$.
                    But since $\xi<\chi_f$ we would have
                    $s_f(\geq\xi)\geq\varepsilon_f$ and hence
                    $s_f(\geq\xi)=\varepsilon_f$, and therefore
                    $\varepsilon_f$ would be singular. Also, we would have $\zeta<s_f(\geq\xi)=\varepsilon_f$ for
                    all $\zeta$ in the support of $s_f$, so $s_f(\varepsilon_f)=0$, in contradiction with ($\varepsilon$reg).
                    So we must have
                    $|h[Y^f_{\geq\xi}]|=s_f(\geq\xi)$, and since
                    $|(h\circ f)\inv[y]|\geq \xi$ for all $y\in
                    h[Y^f_{\geq\xi}]$ we are done.
                \end{proof}

                \begin{lem}\label{LEM:MIDSUPP:sf(>=xi)StableConditions}
                    Let $f\in\Oo$ satisfy ($\varepsilon$), (scont), and ($\varepsilon$reg), and let
                    $1\leq \xi<\chi_f$. Then $s_g(\geq\xi)= s_f(\geq\xi)$
                    for all $g\in\fwithS\sm\S$.
                \end{lem}
                \begin{proof}
                    If $\ef=1$ then there is nothing to show, so we
                    may assume that $\ef$ is infinite, by ($\varepsilon$). Then
                    $s_f(\geq\xi)$ is infinite for all $\xi<\chi_f$.
                    Now if $\xi$ is infinite and regular, or if $\xi\leq 2$, then the assertion is a direct consequence
                    of Lemmas \ref{LEM:UBOUND:sg(>=xi)} and \ref{LEM:MIDSUPP:LBOUND:sf(>=xi)}.
                    If $\xi$ is singular or finite and greater than 2,
                    then there exists $\zeta<\xi$ infinite and regular or
                    equal to 2 such that $s_f(\geq\zeta)=s_f(\geq\xi)$, by (scont).
                    We have
                    $s_g(\geq\xi)\leq
                    s_g(\geq\zeta)=s_f(\geq\zeta)=s_f(\geq\xi)$, and $s_g(\geq \xi)\geq s_f(\geq\xi)$ by Lemma
                    \ref{LEM:MIDSUPP:LBOUND:sf(>=xi)}.
                \end{proof}

                \begin{lem}\label{LEM:CHI:chiIsRegular}
                    If $f\in\Oo$ satisfies (scont), then $\chi_f$ is
                    infinite and regular, or $\chi_f\leq 2$.
                \end{lem}
                \begin{proof}
                    If $\chi_f$ was singular or finite and greater than two, then (scont) would imply
                    that there exists $\zeta<\chi_f$ such that
                    $s_f(\geq\zeta)=s_f(\geq\chi_f)$, contradicting that $\chi_f$ is the minimal cardinal $\zeta\leq\kappa$
                    such that $s_f(\geq\zeta)\leq\lambda$ for some
                    $\lambda\in\supp(s_f)$.
                \end{proof}

                \begin{lem}\label{LEM:CHI:ConditionsForChiStable}
                    Let $f\in\Oo$ satisfy ($\varepsilon$), (scont), and ($\varepsilon$reg).
                    Then $\chi_g=\chi_f$ for all $g\in\fwithS\sm\S$.
                \end{lem}
                \begin{proof}
                    Using ($\varepsilon$), we assume that $\ef$ is infinite. By (scont) and Lemma
                    \ref{LEM:CHI:chiIsRegular}, $\chi_f\leq 2$ or $\chi_f$ is
                    infinite and regular. Assume $\chi_g<\chi_f$. By Lemma
                    \ref{LEM:MIDSUPP:LBOUND:sf(>=xi)},
                    $s_g(\geq\chi_g)\geq s_f(\geq\chi_g)$. But
                    $s_f(\geq\chi_g)>\lambda$ for all
                    $\lambda\in\supp(s_f)$, and therefore also
                    for all $\lambda\in\supp(s_g)$, since $\varepsilon_g=\ef$ by ($\varepsilon$)
                    and Lemma \ref{LEM:EPSILON:eps=1OrInfinite->Stable}, and
                    since ($\varepsilon$reg) and Lemma \ref{LEM:UBOUND:regular} in addition imply that $s_g(\ef)>0$ only if
                    $s_f(\ef)>0$. Thus,
                    $s_g(\geq\chi_g)>\lambda$ for all
                    $\lambda\in\supp(s_g)$, contradicting the definition
                    of $\chi_g$. Assume now that $\chi_g>\chi_f$. Then
                    $s_g(\geq\chi_f)>\lambda$ for all $\lambda\in\supp(s_g)$,
                    and hence also for all $\lambda\in\supp(s_f)$. In particular,
                    $s_g(\geq\chi_f)\geq\ef$ is infinite; thus by Lemma \ref{LEM:UBOUND:sg(>=xi)} we have
                    that $s_f(\geq\chi_f)$ is infinite as well and
                    $s_g(\geq\chi_f)\leq s_f(\geq\chi_f)$, so $s_f(\geq\chi_f)>\lambda$ for
                    all $\lambda\in\supp(s_f)$, in contradiction with the
                    definition of $\chi_f$.
                \end{proof}

                \begin{lem}\label{LEM:MIDSUPP:standardizeSf}
                    Let $f\in\Oo$ satisfy ($\mu$), ($\sigma$), (n), ($\varepsilon$),
                    ($\varepsilon$reg), and assume $\ef\leq\nu_f$.
                    There exists $g\in\fwithS$ such that
                    $s_g(\xi)=s_g(\geq\xi)=s_f(\geq\xi)$ for all $\varepsilon_f'< \xi<\chi_f$ and such that $s_g(\xi)=s_f(\xi)$ for
                    all $\xi\leq\ef'$ and all $\xi\geq\chi_f$.
                \end{lem}

                \begin{proof}
                    We can assume that $\varepsilon_f$ is infinite;
                    for otherwise, $\supp(s_f)=\{1\}$ by condition
                    ($\varepsilon$), and the lemma would be trivial. Also, we assume $\efp<\chi_f$, so in particular
                    $\chi_f>\mu_f$. Define $\delta'\leq \kappa$ to be
                    minimal with the property that $s_f(\zeta)<\ef$
                    for all $\zeta\geq\delta'$, if such a cardinal exists,
                    and to equal $\chi_f$ otherwise. Set $\delta=\min\{\delta', \chi_f\}$.
                    Define
                    $H_\xi=\{\xi\leq\zeta<\delta:s_f(\zeta)\geq\ef\}$, for all
                    $\xi\in(\efp,\delta)$. We claim that $\sum_{\zeta\in
                    H_\xi}s_f(\zeta)=s_f(\geq\xi)$.
                    Indeed, if $s_f(\geq\xi)=\ef$, then there exists $\xi\leq\zeta<\delta$ with $s_f(\zeta)=\ef$, since
                    $\xi<\delta$, so our claim is true. If $s_f(\geq\xi)>\ef$ and our claim did not
                    hold, then we would have $\sum_{\zeta\in
                    [\xi,\ef]\sm
                    H_\xi}s_f(\zeta)=s_f(\geq\xi)$, since $s_f(\geq\xi)\geq\ef$ is infinite.
                    But then
                    $$
                        s_f(\geq\xi)=\sum_{\zeta\in
                        [\xi,\ef]\sm H_\xi}s_f(\zeta)\leq \sum_{\zeta\in
                        [\xi,\ef]\sm H_\xi}\ef\leq\ef,
                    $$
                    contradicting $s_f(\geq\xi)>\ef$.
                    For every $\varepsilon_f'<\xi<\delta$ with $s_f(\xi)\geq \ef$,
                    write $Y^f_\xi$
                    as a disjoint union $\bigcup_{\zeta\leq \xi} Y_{\xi,\zeta}$ in such a way that
                    $|Y_{\xi,\zeta}|=|Y_{\xi}^f|=s_f(\xi)$. This is possible
                    as $s_f(\xi)\geq\varepsilon_f$ is infinite. Set
                    $Y'_\zeta=\bigcup_{\xi\in H_\zeta}Y_{\xi,\zeta}$, for
                    all $\zeta\in (\efp,\delta)$. Then $Y'_\zeta\sub Y^f_{\geq\zeta}$ and $|Y'_\zeta|=\sum_{\xi\in
                    H_\zeta}|Y_{\xi,\zeta}|=\sum_{\xi\in
                    H_\zeta}s_f(\xi)=s_f(\geq\zeta)$.\\
                    If $\delta<\chi_f$, then $s_f(\geq\delta)\geq\ef$ by definition of $\chi_f$, and
                    $s_f(\zeta)<\ef$ for all $\zeta\geq\delta$ implies $s_f(\geq\delta)\leq\ef$, so $s_f(\geq\delta)=\ef$.
                    In that case we must have $\chi_f=\varepsilon_f$: Indeed,
                    $\delta<\chi_f$ implies $s_f(\ef)=0$, so $\ef$ must be regular by
                    ($\varepsilon$reg). Now observe that
                    $\ef=s_f(\geq\delta)=\sum_{\delta\leq
                    \zeta}s_f(\zeta)=\sum_{\delta\leq
                    \zeta<\chi_f}s_f(\zeta)$, which is only possible if
                    $\chi_f=\ef$ by the regularity of $\ef$.
                    Because $\chi_f=\ef$ is a limit cardinal, the support of $s_f$ above
                    $\delta$ is unbounded in $\ef$, and we can find disjoint sets
                    $Y'_\xi\sub Y^f_{\geq\xi}$, for all $\delta\leq\xi<\chi_f$, with $|Y'_\xi|=\ef$.\\
                    Write $Y$ for the union over all $Y'_\xi$ with
                    $\varepsilon_f'<\xi<\chi_f$.
                    Now fix for all $\efp<\xi<\chi_f$ and all $y\in
                    Y_\xi'$ a set $B_{y}\sub f\inv[y]$ with
                    $|B_y|=\xi$,
                    and set $C_y= f\inv [y]\sm B_y$.  Write $B=\bigcup_{y\in
                    Y}B_y$, and $C=\bigcup_{y\in Y} C_y$. Fix $S\sub
                    Y_{\mu_f}^f$ large and such that $Y_{\mu_f}^f\sm S$ is still
                    large, and let $\alpha$ map
                    $S$ bijectively
                    onto $B\cup f\inv[S]$. Let $\alpha$
                    map all $y\in Y^f_\xi$, where $\mu_f<\xi\leq\efp$ or $\xi\geq \chi_f$,
                    into $f\inv[y]$. Write $F=Y^f_{>1,<\aleph_0}\cap (Y^f_{\leq\efp}\cup Y^f_{\geq\chi_f})$, and let $F^*$ consist of
                    those elements of $f\inv[F]$ which are not yet in the range of $\alpha$.
                    Set $D=Y^f_{>\efp,<\chi_f}\sm Y$, and $D^*=f\inv[D]$. Let $\alpha$ map a suitable part of
                    $X\sm f[X]$ bijectively
                    onto $C\cup D^*\cup F^*$. We can do that since $|C\cup D^*\cup F^*|\leq
                    |f\inv[Y^f_{>\efp}]\cup f\inv[Y^f_{>1,<\aleph_0}]|\leq\nu_f+\nu_f=\nu_f$,
                    by the definition of $\efp$ and by (n).
                    Choose moreover $T\sub Y^f_{\mu_f}\sm S$ with $|T|=|Y^f_{>\efp,<\chi_f}|$ and so that
                    $Y^f_{\mu_f}\sm(S\cup T)$
                    is still large, and let $\alpha$ map all
                    $y\in Y^f_{>\efp,<\chi_f}$ into $f\inv[T]$ in
                    such a way that every kernel class of $f\inv[T]$ is hit exactly once.
                    Extend $\alpha$
                    to a bijection, and set $g=f\circ\alpha\circ f$; we
                    can do that since $\alpha$ is not defined on $Y^f_{\mu_f}\sm(S\cup
                    T)$ and its range is disjoint from $f\inv[Y^f_{\mu_f}\sm(S\cup
                    T)]$.
                    We calculate $|g\inv[y]|$ for all $y\in f[X]$. Assume
                    first that $y\in Y$, and say that $y\in
                    Y_{\xi}'$, where $\varepsilon_f'<\xi<\chi_f$.
                    Then $|g\inv [y]|=|(\alpha\circ f)\inv
                    [B_y]|=\mu_f\cdot |B_y|=|B_y|=\xi$,
                    since $\mu_f$ is one or infinite, and
                    $\xi\geq\mu_f$. Assume
                    now that $y\nin Y$. If $y\in Y^f_\xi$
                    for some infinite $\mu_f<\xi\leq\efp$ or $\xi\geq \chi_f$, then
                    $|g\inv[y]|=|f\inv[y]|=\xi$, since $\alpha(y)\in f\inv[y]$ and since $|(\alpha\circ f)\inv[z]|\in\{0,\mu_f,\xi\}$ for
                    all $z\in f\inv[y]$. If $y\in Y^f_\xi$
                    for some finite $\mu_f<\xi\leq\efp$ or
                    $\chi_f\leq\xi\leq\ef$, then $|g\inv[y]|=|f\inv[y]|=\xi$, since
                    $\alpha(y)\in f\inv[y]$ and since $|(\alpha\circ f)\inv[z]|=0$ for all $z\in f\inv[y]$ except
                    $\alpha(y)$.
                    If $y\in Y^f_\xi$ for
                    some $\efp<\xi<\chi_f$ but $y\nin Y$, then $y\in D$
                    and therefore
                    $|g\inv[y]|=0$. If $y\in S$, then $|g\inv[y]|=\mu_f\cdot\mu_f=\mu_f$. If $y\in T$, then
                    there exists exactly one
                    $z\in f\inv[y]\cap \alpha[Y]$, and $|(\alpha\circ f)\inv[w]|\in\{0,\mu_f\}$ for all $w\in
                    f\inv[y]$ except that $z$, so we have $\ef<|g\inv[y]|<\chi_f$. If $y\in Y^f_{\mu_f}\sm (S\cup T)$ then
                    $|g\inv [y]|\in\{0,\mu_f\}$. Therefore we have that
                    for all $\mu_f<\xi\leq\efp$ and
                    all $\xi\geq \chi_f$, $Y^g_\xi=Y^f_\xi$ and thus
                    $s_g(\xi)=s_f(\xi)$; $s_g(\mu_f)=|S|=\kappa=s_f(\mu_f)$; and
                    finally, for all $\efp<\xi<\chi_f$ we have
                    $|g\inv[y]|=\xi$ iff $y\in Y'_\xi$ or $y\in
                    f\circ\alpha[Y^f_\xi]\sub T$, so
                    $s_g(\xi)=|Y'_\xi|+|Y^f_\xi|=s_f(\geq\xi)$; this also
                    implies $s_f(\geq\xi)=s_g(\xi)\leq s_g(\geq\xi)=\sum_{\xi\leq\zeta\leq\ef}s_g(\zeta)\leq
                    \sum_{\xi\leq\zeta\leq\ef}s_f(\geq \zeta)\leq\ef\cdot
                    s_f(\geq\xi)=s_f(\geq\xi)$, so
                    $s_g(\geq\xi)=s_g(\xi)=s_f(\geq\xi)$.
                \end{proof}

                \begin{lem}\label{LEM:MIDSUPP:modifyFinitePartOfSf}
                    Let $f\in\Oo$ satisfy ($\mu$), ($\sigma$), (n), ($\varepsilon$), (scont), ($\varepsilon$reg),
                    and $\ef\leq\nu_f$. Let $p\in\Oo$ be so that
                    $\chi_p=\chi_f$ and
                    $s_p(\geq\xi)=s_f(\geq\xi)$ for all $\efp<\xi<\chi_f$.
                    Then there exists $g\in\fwithS$ such that
                    $s_g(n)=s_p(n)$ for all finite $\efp<n<\chi_f$, such that
                    $s_g(\xi)=s_g(\geq\xi)=s_f(\geq \xi)$ for all infinite
                    $\efp<\xi<\chi_f$, and such that $s_g(\xi)=s_f(\xi)$ for
                    all other $\xi\leq\kappa$.
                \end{lem}
                \begin{proof}
                    We assume that $\ef$ is infinite, using ($\varepsilon$); hence, $\nu_f\geq\ef$ is infinite, too.
                    By Lemma
                    \ref{LEM:MIDSUPP:standardizeSf}, we may
                    assume that $s_f(\xi)=s_f(\geq\xi)$ for all
                    $\varepsilon_f'<\xi<\chi_f$, since this modification
                    obviously does not change the conditions $f$ satisfies,
                    nor the values of $\ef$, $\efp$, $\chi_f$, and $\nu_f$, the latter one staying unchanged by
                    Lemma \ref{LEM:NU:nu_g}. Then there is nothing left to show if
                    $\efp$ is infinite, so we assume it is finite and therefore
                    $\efp=\mu_f=1$ by (n) and ($\mu$). Also, we can assume $\chi_f>2$,
                    so $\chi_f$ is infinite by (scont) and
                    Lemma \ref{LEM:CHI:chiIsRegular}.\\
                    Because $\chi_f>2$, we have that $s_f(\geq 2)\geq\ef$ is
                    infinite. Fix for every $1<n<\aleph_0$ a set
                    $Z_n\subseteq Y^f_n$ with $|Z_n|=s_p(n)$, and set
                    $W_n=Y^f_n\sm Z_n$. Set $Z=\bigcup_{1<n<\aleph_0} Z_n$
                    and $W=\bigcup_{1<n<\aleph_0} W_n$.\\
                    Assume first that $s_f(\geq
                    2)=s_f(\geq\aleph_0)$; then $\aleph_0<\chi_f$ and hence $s_f(\geq
                    2)=s_f(\aleph_0)$. Let $\alpha$ map
                    $y$ into $f\inv[y]$, for all $y\in Z$ and all $y\in
                    Y^f_{\geq\aleph_0}$. Now let $\alpha$
                    map $W$ into $f\inv[Y^f_{\aleph_0}]$ in such a way that it
                    stays injective; since $|W|\leq s_p(\geq
                    2)=s_f(\aleph_0)$, there is enough room to do so.
                    Let $\alpha$ map a suitable part of $X\sm
                    f[X]$ bijectively onto the set of those elements of
                    $f\inv[Z\cup W]$ which are not yet in the image of
                    $\alpha$. The set $X\sm f[X]$ is large enough as
                    $|f\inv[Z\cup W]|\leq |f\inv [Y^f_{>\efp}]|\leq\nu_f$. Choose
                    $S\sub Y^f_1$ large and such that $Y^f_1\sm S$ is large, and map $S$ bijectively onto
                    $f\inv[S]$. Extend the partial injection $\alpha$ to a
                    bijection and set $g=f\circ\alpha\circ f$. If $y\in Y^f_{1}$, then $|g\inv[y]|\leq
                    1$; if $y\in S$, then $|g\inv[y]|= 1$.
                    If $y\in
                    Z_n$, then $|g\inv [y]|=|f\inv[y]|=n$. If $y\in W$,
                    then $f\inv[y]\sub \alpha[X\sm f[X]]$, so
                    $|g\inv[y]|=0$. If $y\in Y^f_{\xi}$ for an infinite $\xi$, then
                    $|g\inv[y]|=\xi$. Therefore,
                    $s_g(\xi)=s_f(\xi)$ for all infinite $\xi\leq\kappa$,
                    $s_g(n)=|Z_n|=s_p(n)$ for all $1<n<\aleph_0$, and
                    $s_g(1)=s_f(1)=\kappa$ and we are done.\\
                    So assume now
                    that $s_f(\geq 2)>s_f(\geq\aleph_0)$.
                    Then
                    also $s_p(\geq 2)>s_p(\geq\aleph_0)$, by the assumptions $s_p(\geq\xi)=s_f(\geq\xi)$ for all
                    $1<\xi<\chi_f$ and $\chi_p=\chi_f$. Therefore, $|Z|=s_p(\geq
                    2,<\aleph_0)=s_p(\geq 2)=s_f(\geq 2)=s_f(\geq
                    2,<\aleph_0)$, and we can find a
                    bijection $\gamma$ from $Y^f_{>1,<\aleph_0}$ onto $Z$ such that
                    whenever $z\in Y^f_n$, then $\gamma(z)\in Z_{j}$ for
                    some $j\geq n$. For every such $z$, we fix a set
                    $B_z\subseteq f\inv[\gamma(z)]$ such that $|B_z|=j-n$, and an
                    element $b_z\in f\inv[\gamma(z)]\sm B_z$. Set
                    $B=\bigcup\{B_z:z\in Y^f_{>1,<\aleph_0}\}$.
                    Let $\alpha$ map every
                    $z\in Y^f_{>1,<\aleph_0}$ to $b_z$. Fix a large set $S\sub
                    Y^f_1$ such that $Y^f_1\sm S$ is large, and let $\alpha$ map $S$ bijectively onto
                    $f\inv[S]\cup B$. Now let $\alpha$ map a suitable part of $X\sm
                    f[X]$ onto the set of those elements of
                    $f\inv[Y^f_{>1,<\aleph_0}]$ which are not in the
                    range of $\alpha$; this is possible as $f\inv[Y^f_{>1,<\aleph_0}]\sub f\inv[Y^f_{>\efp}]$ is not larger
                    than $X\sm f[X]$. Map all
                    $y\in Y^f_{\geq\aleph_0}$ into $f\inv[y]$. Extend $\alpha$ to a
                    bijection and set $g=f\circ\alpha\circ f$. Now if
                    $y\in Y^f_1$, then $|g\inv[y]|\leq 1$, and $|g\inv[y]|= 1$ for all $y\in S$. If $y\in Z_j$ for some
                    $1<j<\aleph_0$,
                    then there exist $1<n\leq j$ and $z\in Y^f_n$ with
                    $\gamma(z)=y$, and we have
                    $|g\inv[y]|=|(\alpha\circ f)\inv[b_z]|+|(\alpha\circ f)\inv[B_z]|=|f\inv[z]|+ 1\cdot |B_z|=n+(j-n)=j$. If
                    $y\in Y^f_{>1,<\aleph_0}\sm Z$, then $|g\inv[y]|=0$.
                    If $y\in Y^f_\xi$ for some infinite $\xi\leq\kappa$,
                    then $|g\inv[y]|=|f\inv[y]|=\xi$. Therefore,
                    $s_g(j)=|Z_j|=s_p(j)$ for all $1<j<\aleph_0$, $s_g(\xi)=s_f(\xi)$ for
                    all infinite $\xi$, and $s_g(1)=s_f(1)=\kappa$.
                \end{proof}

                \begin{lem}\label{LEM:MIDSUPP:modifySfArbitrarily}
                    Let $f\in\Oo$ satisfy ($\mu$), ($\sigma$), (n), ($\varepsilon$), (scont), ($\varepsilon$reg),
                    and $\ef\leq\nu_f$. Let $p\in\Oo$ be so
                    that $\chi_p=\chi_f$ and
                    $s_p(\geq\xi)=s_f(\geq\xi)$ for all $\efp<\xi<\chi_f$.
                    Then there exists $g\in\fwithS$ such that
                    $s_g(\xi)=s_p(\xi)$ for all $\varepsilon_f'<\xi
                    <\chi_f$, and $s_g(\xi)=s_f(\xi)$ for all
                    $\xi\leq\varepsilon_f'$ and all $\xi\geq\chi_f$.
                \end{lem}
                \begin{proof}
                    We can assume that $\ef$ is infinite, for otherwise there is nothing to show; hence, $\nu_f\geq\ef$ is
                    infinite, too. We also assume that $\efp<\chi_f$, so in particular $\chi_f>\mu_f$. By Lemma \ref{LEM:MIDSUPP:modifyFinitePartOfSf}, we may
                    assume that $s_f(n)=s_p(n)$ for all finite $\efp<n<\chi_f$, and that
                    $s_f(\xi)=s_f(\geq \xi)$ for all infinite
                    $\efp<\xi<\chi_f$.
                    Fix for every infinite $\efp<\xi<\chi_f$ a set $Z_\xi\sub
                    Y^f_\xi$ with $|Z_\xi|=s_p(\xi)$. If $s_p(\xi)=s_f(\xi)=|Y^f_\xi|$, then
                    we shall have $Z_\xi=Y^f_\xi$. For every $z\in
                    Z_\xi$, write $f\inv[z]=\{z_\zeta:\zeta<\xi\}$ (observe that here, the $\zeta<\xi$ are of course
                    all ordinals below
                    $\xi$ and not only cardinals). Write
                    $A_\xi^\zeta=\{z_\zeta: z\in Z_{\xi}\}\sub f\inv[Z_\xi]$, for all
                    $\zeta<\xi$ (so the index $\xi$ is a cardinal, and the index $\zeta$ an ordinal below
                    $\xi$). Then $|A_\xi^\zeta|=|Z_\xi|$.
                    We define a
                    partial injection $\alpha$ in the following way: Map each $Z_\xi$ bijectively onto
                    $A_\xi^0$. Next we define $\alpha$ on $Y_\xi\sm
                    Z_\xi$; this is only necessary if $Y_\xi\sm
                    Z_\xi\neq\emptyset$, which happens if and only if $s_p(\xi)<s_f(\xi)$.
                    In that case, set
                    $Z_{>\xi}=\bigcup_{\zeta>\xi} Z_\zeta$. Because
                    $s_p(\xi)<s_f(\xi)$ but $s_p(\geq\xi)=s_f(\geq\xi)=s_f(\xi)$, we have $s_p(>\xi)=s_f(\xi)$;
                    in particular, this implies $s_p(>\xi)=s_p(\geq\xi)>s_p(\geq\chi_p)$,
                    and together with $\chi_p=\chi_f$ we infer $s_p(>\xi,<\chi_f)=s_p(>\xi)$. Hence,
                    $|Z_{>\xi}|=s_p(>\xi,<\chi_f)=s_p(>\xi)=s_f(\xi)$. Therefore,
                    $|\bigcup_{\zeta>\xi} A_\zeta^\xi|=|Z_{>\xi}|=s_f(\xi)=|Y^f_\xi\sm
                    Z_\xi|$. Let $\alpha$ map $Y^f_\xi\sm
                    Z_\xi$ bijectively onto $\bigcup_{\zeta>\xi}
                    A_\zeta^\xi$. The function $\alpha$ is injective.
                    Indeed, it is injective on each $Z_\xi$ and $Y^f_\xi\sm
                    Z_\xi$ by definition, and $f[Z_\xi]=A_\xi^0$, and $f[Y^f_\xi\sm
                    Z_\xi]=\bigcup_{\zeta>\xi}
                    A_\zeta^\xi$, so the ranges of those
                    injective parts are disjoint. Extend $\alpha$ by
                    mapping $y$ into $f\inv[y]$ for all $y\in Y^f_\xi$
                    for which $\mu_f<\xi\leq\varepsilon'_f$ or
                    $\xi\geq\chi_f$, or for which $\efp<\xi<\chi_f$ is finite. Extend $\alpha$ by
                    mapping a suitable part of $X\sm f[X]$ bijectively onto those
                    elements of $f\inv [Y^f_{>\varepsilon_f'}\cup Y^f_{>1,<\aleph_0}]$ which are
                    not yet in the range of $\alpha$; this can be done as $f\inv [Y^f_{>\varepsilon_f'}]$ and
                    $f\inv[Y^f_{>1,<\aleph_0}]$ are not larger than $\nu_f$, the first one by
                    definition of $\efp$ and the second
                    one by condition (n). Choose $S\sub
                    Y^f_{\mu_f}$ large and such that $Y^f_{\mu_f}\sm S$ is
                    large, and let $\alpha$ map $S$ bijectively onto
                    $f\inv[S]$. Extend $\alpha$ to a
                    bijection and set $g=f\circ\alpha\circ f$.\\
                    We calculate $|g\inv[z]|$ for all $z\in X$.
                    If $z\in Z_\xi$ for some infinite $\efp<\xi<\chi_f$,
                    then $\alpha(z)\in f\inv[z]$, so $g\inv[z]\supseteq
                    f\inv[z]$ and hence $|g\inv[z]|\geq\xi$. On the other
                    hand, if $\alpha(w)=z_\zeta\in f\inv[Z_\xi]$, then
                    $\alpha(w)\in A_\xi^\zeta$ by definition of that set.
                    If $\zeta=0$, then $w\in Z_\xi$, so $|f\inv
                    [w]|=\xi$, and if $\xi>\zeta>0$, then
                    $w\in Y^f_\zeta\sm Z_\zeta$, so $|f\inv
                    [w]|<\xi$. Thus
                    $|g\inv[z]|\leq\xi$ and whence, $|g\inv[z]|=\xi$.\\
                    If $z\in Y^f_\xi\sm Z_\xi$ for
                    some infinite $\varepsilon_f'<\xi<\chi_f$, then $f\inv[z]\cap (\alpha\circ
                    f[X])=\emptyset$, and so $g\inv[z]=\emptyset$.
                    If $z\in Y_\xi^f$ for some
                    $\xi\geq\chi_f$, then $(f\circ\alpha)\inv[z]\cap
                    f[X]=\{z\}$, by definition of $\alpha$, and so
                    $|g\inv[z]|=|f\inv[z]|=\xi$. If $z\in Y_\xi^f$ for
                    some infinite
                    $\mu_f<\xi\leq\varepsilon'_f$, then $z\in(f\circ\alpha)\inv[z]$, by definition of $\alpha$,
                    and so
                    $|g\inv[z]|\geq |f\inv[z]|=\xi$. Moreover, if $w\in
                    (f\circ\alpha)\inv[z]$ is distinct from $z$, then $w\in X\sm f[X]$ or $w\in
                    Y^f_{\mu_f}$, and therefore $|g\inv[z]|=\xi$. If $z\in
                    Y^f_\xi$ for some $1<\xi<\aleph_0$, then
                    $(f\circ\alpha)\inv[z]\cap f[X]=\{z\}$, so
                    $|g\inv[z]|=\xi$.
                    If $z\in Y^f_{\mu_f}$, then
                    $|(\alpha\circ f)\inv [w]|\leq\mu_f$ for all $w\in
                    f\inv[z]$, so $|g\inv [z]|\in\{0,\mu_f\}$; if $z\in S$, then $|g\inv [z]|=\mu_f$. Therefore,
                    $s_g(\xi)=s_f(\xi)$ for all $\xi\leq\varepsilon_f'$
                    and all $\xi\geq\chi_f$, and for all finite $\xi$. Also, we have seen that
                    $|g\inv[z]|=\xi$ with $\efp<\xi<\chi_f$ infinite if and only
                    if $z\in Z_\xi$, and thus $s_g(\xi)=s_p(\xi)$.
                \end{proof}

                \begin{lem}\label{LEM:MIDSUPP:MODIFY:arbitraryBelowChi}
                    Let $f\in\Oo$ satisfy  ($\mu$), ($\nu$), ($\sigma$), (s'dec), (n), ($\varepsilon$),
                    (scont), ($\varepsilon$reg), and $\ef\leq\nu_f$.
                    Let $p\in\Oo$ be so
                    that $s_p'=s_f'$, $\mu_p=\mu_f$, $\sigma_p=\sigma_f=\kappa$, and $\nu_p=\nu_f$.
                    Assume moreover that $\chi_p=\chi_f$ and
                    $s_p(\geq\xi)=s_f(\geq\xi)$ for all $\efp<\xi<\chi_f$.
                    Then there exists $g\in\fwithS$ such that
                    $s_g\rest_{[1,\chi_f)}=s_p\rest_{[1,\chi_f)}$ and such that
                    $s_g\rest_{[\chi_f,\kappa]}=s_f\rest_{[\chi_f,\kappa]}$.
                \end{lem}
                \begin{proof}
                    Observe first that (s'inf) holds for $f$ by
                    Lemma \ref{LEM:SUF:(s'inf)} and since $f$
                    satisfies ($\nu$), (s'dec), and $\ef\leq\nu_f$.
                    We have $\efp=\lambda_f'$ since $\ef\leq\nu_f$. By Lemma
                    \ref{LEM:WSUPP:BELOWLAMBDA:addingOfClasses},
                    we can find $h\in\fwithS$ such that $s_h\rest_{[1,\efp)}=s_p\rest_{[1,\efp)}$
                    and such that $s_h\rest_{[\efp,\kappa]}=s_f\rest_{[\efp,\kappa]}$. Observe
                    that $s_h(\efp)=s_f(\efp)=s_p(\efp)$ since $s_p'=s_f'$ and $\sigma_p=\sigma_f$ and
                    since $\efp$ either is an element of $\supp'(s_f)$ or equals
                    $\mu_f$. Clearly, $\chi_h=\chi_f$.
                    Also, because of ($\nu$) and Lemma \ref{LEM:NU:(nu)->NuStable} we have
                    $\nu_h=\nu_f$, and hence $\varepsilon_h'=\efp$. Moreover,
                    $s_h(\geq\xi)=s_f(\geq\xi)$ for all $\xi>\efp$. It is
                    easy to see that $h$ still satisfies the conditions of Lemma
                    \ref{LEM:MIDSUPP:modifySfArbitrarily}. Hence, there is
                    $g\in\cl{\{h\}\cup\S}$ such that $s_g(\xi)=s_p(\xi)$
                    for all $\varepsilon_f'<\xi
                    <\chi_f$, and $s_g(\xi)=s_h(\xi)$ for all other $\xi\leq\kappa$.
                \end{proof}


            \subsubsection{The lion-tail}

                We modify functions $f\in\Oo$ beyond $\chi_f$, thereby completing the case
                $\ef\leq\nu_f$.

                \begin{lem}\label{LEM:WWSUPP:makingSfZero}
                    Let $f\in\Oo$ satisfy ($\mu$), ($\sigma$), (n), ($\varepsilon$), ($\chi$),
                    $\varepsilon_f\leq\nu_f$, and
                    $s_f(\geq\chi_f)>0$. There exists $g\in\Oo$ such that
                    $s_g(\xi)=0$ for all $\chi_f\leq\xi<\varepsilon_f$,
                    such that $s_g(\varepsilon_f)=1$, and such that
                    $s_g(\xi)=s_f(\xi)$ for all other $\xi\leq\kappa$.
                \end{lem}
                \begin{proof}
                    Observe that $\mu_f<\chi_f$, for otherwise
                    $\kappa=s_f(\mu_f)=s_f(\geq\chi_f)$ by ($\sigma$), contradicting ($\chi$).
                    We can assume $\ef>1$ and therefore that $\ef$ is infinite by ($\varepsilon$).
                    Because $s_f(\geq\chi_f)>0$ and ($\chi$) we have that
                    $s_f(\varepsilon_f)>0$. Fix $y\in
                    Y^f_{\varepsilon_f}$, and let $\alpha$ map
                    $Y^f_{\geq\chi_f}$ injectively into $f\inv[y]$.
                    Now let $\alpha$ map every $z\in Y^f_\xi$, where
                    $\mu_f<\xi<\chi_f$, into $f\inv[z]$. Extend $\alpha$ by mapping a suitable part of
                    $X\sm f[X]$ injectively onto the set of those elements of
                    $f\inv[Y^f_{\geq\chi_f}]\cup f\inv[Y^f_{>1,<\aleph_0}]$ which are not yet in the range of $\alpha$;
                    we can do that since $\ef\leq\nu_f$ implies $|f\inv[Y^f_{\geq\chi_f}]|\leq\nu_f$ and (n) implies
                    $|f\inv[Y^f_{>1,<\aleph_0}]|\leq\nu_f$.
                    Choose $S\sub Y^f_{\mu_f}$ large and such that $Y^f_{\mu_f}\sm S$
                    is large and let $\alpha$ map $S$ bijectively onto $f\inv[S]$. Extending
                    $\alpha$ to a bijection, we claim that
                    $g=f\circ\alpha\circ f$ has the desired properties and
                    calculate $|g\inv [z]|$ for all $z\in X$. If $z\in
                    Y^f_{\mu_f}$, then $|g\inv [z]|\in\{0,\mu_f\}$ by
                    construction of $\alpha$ and since $\mu$ is either 1
                    or infinite by ($\mu$); if $z\in S$ then $|g\inv [z]|=\mu_f$. If $z\in
                    Y^f_\xi$, where $\mu_f<\xi<\chi_f$, then $|g\inv
                    [z]|=|f\inv [z]|=\xi$. If $\chi_f\leq \xi\leq \varepsilon_f$ and $z\in Y^f_\xi$,
                    then $|g\inv [z]|=0$ unless $z=y$, in which case $|g\inv
                    [z]|=\ef$. Therefore, if $\chi_f\leq\xi\leq\ef$, then
                    $s_g(\xi)=0$ unless $\xi=\varepsilon_f$, in which case
                    we have $s_g(\xi)=1$. If $\mu_f<\xi<\chi_f$, then
                    $|g\inv[z]|=\xi$ if and only if $z\in Y^f_\xi$, so
                    $s_g(\xi)=s_f(\xi)$. If $\xi>\ef$, then $s_g(\xi)=0$. Finally,
                    $s_g(\mu_f)=s_f(\mu_f)=\kappa$.
                \end{proof}

                \begin{lem}\label{LEM:WWSUPP:makingSfZeroAndSf(eps)=n}
                    Let $f\in\Oo$ satisfy ($\mu$), ($\sigma$), (n),
                    ($\varepsilon$), ($\chi$),
                    $\varepsilon_f\leq\nu_f$, and
                    $s_f(\geq\chi_f)>0$. Let $1\leq n<\aleph_0$. There exists $g\in\Oo$ such that
                    $s_g(\xi)=0$ for all $\chi_f\leq\xi<\varepsilon_f$,
                    such that $s_g(\varepsilon_f)=n$, and such that
                    $s_g(\xi)=s_f(\xi)$ for all other $\xi<\kappa$.
                \end{lem}
                \begin{proof}
                    Using Lemma \ref{LEM:WWSUPP:makingSfZero}, it is enough to
                    show that assuming
                    $s_f(\xi)=0$ for all $\chi_f\leq\xi<\varepsilon_f$,
                    and $s_f(\varepsilon_f)=k$, where $1\leq k<\aleph_0$, we can
                    produce $g\in\fwithS$ such that $s_g(\varepsilon_f)=k+1$ and such that $s_g(\xi)=s_f(\xi)$
                    for all other $\xi\leq\kappa$; this is legitimate
                    as application of Lemma
                    \ref{LEM:WWSUPP:makingSfZero}, as well as
                    increasing the value of $s_f(\ef)$ by $1$,
                    does not change the conditions on $f$ and leaves the values
                    $\chi_f$ and $\ef$ unaltered.
                    To do this, let $\alpha$ map all $z\in
                    Y_\xi^f$, where $\mu_f<\xi<\varepsilon_f$, into $f\inv[z]$. Fix $y \in
                    Y^f_{\mu_f}$ and write
                    $Y^f_{\varepsilon_f}=\{z_1,\ldots,z_k\}$. Let $\alpha$
                    map $z_i$ into $f\inv[z_i]$, $2\leq i\leq k$, and $z_1$ into $f\inv[y]$. Extend $\alpha$ by
                    mapping a suitable part of $X\sm f[X]$ bijectively onto those elements of
                    $f\inv[Y^f_{>1,<\aleph_0}]$ which are
                    not yet in the image of $\alpha$; we can do that by (n). Choose
                    $S\sub Y^f_{\mu_f}$ large and such
                    that $Y^f_{\mu_f}\sm S$ is still large, and let $\alpha$ map
                    $S$ bijectively onto $f\inv[S\cup\{z_1\}]$.
                    Extending $\alpha$ to a bijection, we claim that
                    $g=f\circ\alpha\circ f$ has the desired properties.
                    Indeed, if $\mu_f<\xi<\varepsilon_f$ and $z\in Y^f_\xi$, then
                    $|g\inv[z]|=|f\inv[z]|=\xi$. If $z\in
                    Y^f_{\varepsilon_f}$, then $|g\inv[z]|=|f\inv[z]|=\ef$
                    if $z\in\{z_2,\ldots,z_k\}$, and
                    $|g\inv[z]|=\mu_f\cdot\ef=\ef$ if $z=z_1$ since $f\inv[z_1]\sub\alpha[Y^f_{\mu_f}]$.
                    If $z\in Y^f_{\mu_f}$, then $|g\inv[z]|\in\{0,\mu_f\}$
                    unless $z=y$, in which case
                    $|g\inv[z]|=|f\inv[z_1]|=\ef$; moreover,
                    $|g\inv[z]|=\mu_f$ if $z\in S$. Hence, $s_g(\varepsilon_f)=|\{y,z_1,\ldots,z_k\}|=k+1$,
                    $s_g(\xi)=|Y^f_\xi|=s_f(\xi)$ for
                    all $\mu_f<\xi<\ef$, and
                    $s_g(\mu_f)=|S|=\kappa=s_f(\mu_f)$.
                \end{proof}

                \begin{lem}\label{LEM:WWSUPP:makingSfArbitrary}
                    Let $f\in\Oo$ satisfy ($\mu$), ($\sigma$), (n),
                    ($\varepsilon$),  ($\chi$), $\varepsilon_f\leq\nu_f$, and
                    $s_f(\geq\chi_f)>0$. Let $h\in\Oo$ be so that
                    $\varepsilon_h=\varepsilon_f$ and such that $s_h(\geq\chi_f)$ is finite.
                    Then there exists $g\in\fwithS$ such that
                    $s_g\rest_{[1,\chi_f)}=s_f\rest_{[1,\chi_f)}$ and such that
                    $s_g\rest_{[\chi_f,\kappa]}=s_h\rest_{[\chi_f,\kappa]}$.
                \end{lem}
                \begin{proof}
                    Set $n=s_h(\geq\chi_f)$. By Lemma
                    \ref{LEM:WWSUPP:makingSfZeroAndSf(eps)=n}, we can
                    assume that $s_f(\xi)=0$ for all
                    $\chi_f\leq\xi<\varepsilon_f$, and that
                    $s_f(\varepsilon_f)=n$. Write
                    $Y^f_{\varepsilon_f}=\{z_1,\ldots,z_n\}$ and
                    $Y^h_{\geq\chi_f}=\{y_1,\ldots,y_n\}$. Since $\varepsilon_h=\ef$ we have $Y^h_{\ef}\neq \emptyset$;
                    say without
                    loss of generality that $|h\inv[y_1]|=\varepsilon_f$. For all $\mu_f<\xi<\chi_f$, let $\alpha$ map all
                    $y\in Y_\xi^f$ into $f\inv[y]$.
                    Now let $\alpha$ map $Y^f_{\varepsilon_f}$ injectively into
                    $f\inv[z_1]$. Fix for every $2\leq i\leq n$ a set
                    $Z_i\subseteq f\inv[z_i]$ with $|Z_i|=|h\inv[y_i]|$.
                    Let $\alpha$ map a suitable part of $X\sm f[X]$ bijectively
                    onto the union of $\bigcup_{2\leq i\leq n}f\inv[z_i]\sm Z_i$ with the set of
                    those elements of $f\inv[Y^f_{>1,<\aleph_0}]$ which are not yet in the range of $\alpha$;
                    this can be done as
                    $\ef\leq\nu_f$ and as $f$ satisfies (n). Take $S\sub
                    Y^f_{\mu_f}$ large and so that $Y^f_{\mu_f}\sm S$ is still large, and
                    let $\alpha$ map $S$ bijectively onto $f\inv[S]\cup \bigcup_{2\leq i\leq n}
                    Z_i$. Extend $\alpha$ to a
                    bijection and set $g=f\circ\alpha\circ f$. If
                    $y\in Y_{\mu_f}^f$, then $|g\inv[y]|\in\{0,\mu_f\}$, and if $y\in S$ then $|g\inv[y]|=\mu_f$. If
                    $y\in Y_\xi^f$, where $\mu_f<\xi<\chi_f$, then
                    $|g\inv[y]|=|f\inv[y]|=\xi$. Now assume $y\in
                    Y_{\varepsilon_f}^f$, and say first that $y=z_1$. Then
                    $|g\inv[y]|=\varepsilon_f=|h\inv[y_1]|$. If $y=z_i$,
                    where $2\leq i\leq n$, then $|g\inv [y]|=|Z_i|\cdot\mu_f=|Z_i|=
                    |h\inv
                    [y_i]|$. Thus, we have that $s_g(\xi)=s_f(\xi)$ for
                    all $\xi<\chi_f$, and $s_g(\xi)=s_h(\xi)$ for all
                    $\xi\geq\chi_f$.
                \end{proof}

                \begin{lem}\label{LEM:TAIL:MODIFY:arbitraryEveryWhere}
                    Let $f\in\Oo$ satisfy ($\mu$), ($\nu$), ($\sigma$), (s'dec), (n), ($\varepsilon$),
                    (scont), \cchi, ($\varepsilon$reg), and $\ef\leq\nu_f$ and $s_f(\geq\chi_f)>0$.
                    Let $p\in\Oo$ be so that $s_p'=s_f'$, $\mu_p=\mu_f$, $\sigma_p=\sigma_f=\kappa$, and
                    $\nu_p=\nu_f$. Assume moreover that $\chi_p=\chi_f$, that
                    $s_p(\geq\xi)=s_f(\geq\xi)$ for all $\efp<\xi<\chi_f$, that $\varepsilon_p=\ef$,
                    and that $s_p(\geq\chi_f)>0$ is finite. Then there exists $g\in \fwithS$ such that $s_g=s_p$.
                \end{lem}
                \begin{proof}
                    By Lemma \ref{LEM:MIDSUPP:MODIFY:arbitraryBelowChi}, there exists
                    $q\in\fwithS$ such that
                    $s_q\rest_{[1,\chi_f)}=s_p\rest_{[1,\chi_f)}$ and
                    $s_q\rest_{[\chi_f,\kappa]}=s_f\rest_{[\chi_f,\kappa]}$. This function
                    $q$ obviously still satisfies the conditions of Lemma
                    \ref{LEM:WWSUPP:makingSfArbitrary}; also, $\varepsilon_p=\ef=\varepsilon_q$ and
                    $\chi_p=\chi_f=\chi_q$.
                    Therefore, that lemma implies that $q$
                    together with $\S$ generates a function $g$ such that
                    $s_g\rest_{[1,\chi_f)}=s_q\rest_{[1,\chi_f)}=s_p\rest_{[1,\chi_f)}$
                    and such that
                    $s_g\rest_{[\chi_f,\ef]}=s_p\rest_{[\chi_f,\ef]}$. Hence,
                    $s_g=s_p$.
                \end{proof}

                \begin{prop}\label{PROP:criterionNoClassBeyondNu}
                    Let $f\in\Oo$ be so that $\ef\leq\nu_f$.
                    If $f$ moreover satisfies
                    ($\mu$), ($\nu$), ($\sigma$), ($\rho$), (s'dec), (n), ($\varepsilon$),
                    (scont) and ($\chi$), then it is $\S$-minimal.
                \end{prop}
                \begin{rem}
                    In this situation, $f$ automatically satisfies
                    (\#$\varepsilon$) and
                    ($\lambda$'), as $\ef\leq\nu_f$.
                \end{rem}
                \begin{proof}
                    Let $g\in\fwithS\sm\S$. By Lemmas \ref{LEM:SUF:(epsreg)},
                    \ref{LEM:SUF:(kappa)} and \ref{LEM:SUF:(s'inf)},
                    $f$ satisfies ($\varepsilon$reg), ($\kappa$), and (s'inf). We have $\varepsilon_g=\varepsilon_f$,
                    $\nu_g=\nu_f$, $\mu_g=\mu_f$, $\chi_g=\chi_f$, and $s_g'=s_f'$, by
                    Lemmas \ref{LEM:STABLE:MU},
                    \ref{LEM:CHI:ConditionsForChiStable},
                    and \ref{LEM:STRSUPP:s'dec->stable}, respectively.
                    By Lemma \ref{LEM:MIDSUPP:sf(>=xi)StableConditions},
                    $s_g(\geq\xi)=s_f(\geq\xi)$ for all $\xi<\chi_f$. The latter fact, together with the
                    fact that
                    $\chi_g$ is either infinite and regular or not greater than 2
                    provided by Lemma \ref{LEM:CHI:chiIsRegular}, implies that $g$ satisfies (scont).
                    By Lemma \ref{LEM:STABLE:MU}, $g$ satisfies
                    ($\mu$), ($\nu$), ($\sigma$), ($\varepsilon$), and
                    ($\varepsilon$reg).
                    Because $s_g'=s_f'$, $g$ satisfies (s'dec) and (s'inf), and by
                    Lemma \ref{LEM:STABLE:(Nonu)+(nu)Stable} it
                    satisfies (n).\\
                    Now if $s_f(\geq\chi_f)=0$, then by Lemma \ref{LEM:MIDSUPP:MODIFY:arbitraryBelowChi}
                    we find $h\in\gwithS$ such that $s_h\rest_{[1,\chi_f)}=s_f\rest_{[1,\chi_f)}$ and such that
                    $s_h\rest_{[\chi_f,\kappa]}=s_g\rest_{[\chi_f,\kappa]}$. But by
                    Lemma \ref{LEM:UBOUND:sg(>=xi)}, we have
                    $s_h(\geq\chi_f)=s_g(\geq\chi_f)=s_f(\geq\chi_f)=0$. Hence,
                    $s_h=s_f$ so that since also $\nu_h=\nu_g=\nu_f$ by Lemma \ref{LEM:NU:(nu)->NuStable}, we
                    conclude $f\in\gwithS$.\\
                    If on the other hand $s_f(\geq\chi_f)>0$, then also $s_g(\geq\chi_f)>0$;
                    $s_g(\geq\chi_f)$ is finite by Lemma
                    \ref{LEM:UBOUND:sg(>=xi)}, so $g$
                    satisfies \cchi. Therefore by Lemma
                    \ref{LEM:TAIL:MODIFY:arbitraryEveryWhere} there exists
                    $h\in\gwithS$ such that $s_h=s_f$; since
                    $\nu_h=\nu_g=\nu_f$ we infer $f\in \gwithS$.
                \end{proof}


    \section{Necessities for $\S$-minimality}

        We prove that the conditions of Theorem
        \ref{THM:mainTheorem} are necessary for a function to be
        $\S$-minimal.

        \subsection{Farmers}

            \begin{lem}\label{LEM:MU:fConstantOrLargeRange}
                If $f$ is $\S$-minimal, then it is constant or has large
                range.
            \end{lem}
            \begin{proof}
                If $f$ has small range, then there exists $y\in X$
                with $|f\inv[y]|>|f[X]|$, for $\kappa=\sum_{z\in
                f[X]}|f\inv[z]|$. Let $\alpha$ map $f[X]$ injectively
                into $f\inv [y]$. Since both domain and range of the
                partial function $\alpha$ are co-large, we can extend
                it to a bijection on $X$. The function
                $g=f\circ\alpha\circ f$ is constant and an element of
                $\cl{\{f\}\cup\S}$. Since $f$ is $\S$-minimal, we must have
                $f\in\cl{\{g\}\cup\S}$, which is only possible if $f$
                is constant itself.
            \end{proof}
            \begin{lem}\label{LEM:MU:rhoIsSmall}
                If $f$ is $\S$-minimal, then it satisfies ($\rho$), i.e. $\varrho_f<\kappa$.
            \end{lem}
            \begin{proof}
                Assume to the contrary that $Y^f_{>\mu_f}$ is large. Let
                $Z\subseteq Y^f_{>\mu_f}$ be large and so that $Y^f_{>\mu_f}\sm Z$ is large as well, and let $\alpha$
                map $Y^f_{\mu_f} \cup Z$ bijectively onto $f\inv[Z]$. Both range and
                domain of $\alpha$ are co-large, so we can extend it
                to a permutation on $X$. Now all kernel classes of $g=f\circ\alpha\circ f$
                are strictly larger than $\mu_f$. Indeed, assume
                $|g\inv [y]|=\mu_f$ for some $y\in X$. Then there
                exists $z\in f\inv [y]\cap (\alpha\circ f) [X]$; for
                this $z$ we must have $|(\alpha\circ f)\inv
                [z]|=\mu_f$. By construction of $\alpha$ we conclude
                that $z\in f\inv[Z]$, and so $y\in Z$. But then
                $|g\inv [y]|\geq |f\inv [y]\cap (\alpha\circ f)[X]|=|f\inv
                [y]|>\mu_f$, contradiction.
                So indeed $\mu_g>\mu_f$ and
                we cannot get back $f$ from $g$ and $\S$, contradicting that
                $f$ be $\S$-minimal.
            \end{proof}
            \begin{lem}\label{LEM:MU:sigmaIsLarge}
                If $f$ is $\S$-minimal and nonconstant, then it satisfies ($\sigma$), i.e. $\sigma_f=\kappa$.
            \end{lem}
            \begin{proof}
                We know from Lemma \ref{LEM:MU:fConstantOrLargeRange} that $f$ has large range. Therefore,
                $\sigma_f+\rho_f=|f[X]|=\kappa$. Since $\rho_f<\kappa$ by Lemma \ref{LEM:MU:rhoIsSmall}, we infer
                $\sigma_f=\kappa$.
            \end{proof}

            \begin{lem}\label{LEM:MU:muIsOneOrInfinite}
                If $f$ is $\S$-minimal, then it satisfies ($\mu$), i.e. $\mu_f=1$ or $\mu_f$ is infinite.
            \end{lem}
            \begin{proof}
                Assume that $\mu_f$ is finite but not equal to $1$. Then $f$ is nonconstant, and therefore
                $Y^f_{\mu_f}$ is large
                by Lemma \ref{LEM:MU:sigmaIsLarge}. Let
                $S\subseteq Y^f_{\mu_f}$ be large and such that $Y^f_{\mu_f}\sm S$ is still large and let $\alpha$
                map $S$ bijectively onto $f\inv [S]$. Both domain and
                range of $\alpha$ are co-large, so we can extend it
                to a bijection on $X$. Set $g=f\circ \alpha\circ f$.
                Then for all $y\in S$, $|g\inv[y]|=\mu_f^2>\mu_f$. Thus,
                $s_g(>\mu_f)=\kappa$. Now if $\mu_g>\mu_f$, then obviously $f\nin\gwithS$, so
                $f$ is not $\S$-minimal;
                if on the other hand $\mu_g=\mu_f$, then Lemma \ref{LEM:MU:rhoIsSmall}
                implies that $g$ is not $\S$-minimal as $\rho_g=s_g(>\mu_g)=\kappa$, hence in that case
                $f$ is not $\S$-minimal either, contradiction.
            \end{proof}

        \subsection{The return of the man who wasn't there}

            \begin{lem}\label{LEM:NU:(mu=1&NuInfinite)ORMuInfinite}
                If $f$ is $\S$-minimal, then it satisfies ($\nu$), i.e. if $\mu_f$ is finite, then $\nu_f$ is
                infinite or zero.
            \end{lem}
            \begin{proof}
                If $\mu_f$ is finite, then $\mu_f=1$ by Lemma
                \ref{LEM:MU:muIsOneOrInfinite}. Assume that in
                this situation, $0<\nu_f<\aleph_0$. Fix $y\in Y^f_1$, and
                choose $\alpha\in\S$ so that $f\inv[y]\cap(\alpha\circ
                f[X])=\emptyset$. Set $g=f\circ \alpha\circ f$; then $\nu_g\geq|(X\sm
                f[X])\cup\{y\}|=\nu_f+1$, in contradiction with the
                obvious fact that if $f$ is $\S$-minimal, then $\nu_g=\nu_f$ for all $g\in\fwithS\sm\S$.
            \end{proof}

        \subsection{The revenge of the dwarf-box}

            \begin{lem}\label{LEM:FINSUPP:nuFinite->sf(n)=0,nuInfinite->sf(n)<=nu}
                Let $f$ be $\S$-minimal. If $\nu_f$ is finite,
                then $s_f(n)=0$ for all $1<n<\aleph_0$. If $\nu_f$ is
                infinite, then $s_f(n)\leq\nu_f$ for all
                $1<n<\aleph_0$. In particular, there exist no finite
                cardinals in the strong support of $s_f$ and $f$ satisfies (n).
            \end{lem}

            \begin{proof}
                The lemma is trivial if $\mu_f$ is infinite, so we can assume $\mu_f=1$ by
                Lemma \ref{LEM:MU:muIsOneOrInfinite}. Then Lemma
                \ref{LEM:NU:(mu=1&NuInfinite)ORMuInfinite} implies that
                $\nu_f$ is zero or infinite, since $f$ satisfies \cnu.
                So all we have to show that $s_f(n)\leq\nu_f$ for all $1<n<\aleph_0$.
                Suppose there is $1<n<\aleph_0$ with
                $s_f(n)>\nu_f$, and let $n$ be minimal with this property.
                Choose $Z\sub Y^f_1$ such that $|Z|=|Y^f_n|$ and such that $Y^f_1\sm Z$ is large; we can do this
                since $Y^f_1$ is large by Lemma \ref{LEM:MU:sigmaIsLarge}.
                For every $y\in Y^f_n$, let $\alpha$ map $y$ into
                $f\inv[y]$, and let it map exactly one
                element of $Z$ into $f\inv[y]$ in such a way that it stays injective.
                Extend the mapping to a bijection on
                $X$ and set $g=f\circ\alpha\circ f$. Observe that
                we must have
                $\nu_g=\nu_f$ since $f$ is $\S$-minimal.
                We claim that for
                all $1< k\leq n$, $s_g(k)\leq \nu_f$. Indeed, for
                $1<k<n$ this follows from Lemma \ref{LEM:UBOUND:FINITE:sg},
                since $s_f(i)\leq\nu_f$ for all $1<i\leq k$ and since
                $\nu_f$ is zero or infinite. Now
                assume $|g\inv [y]|=n$ for some $y\in X$.
                If $|f\inv [y]|>n$, then $f\inv[y]\sm (\alpha\circ f[X])$ must be non-empty, which happens at
                most $\nu_f$ times. If $|f\inv [y]|=n$, then $\alpha(y)\in
                f\inv[y]$, so $g\inv[y]$ contains $f\inv[y]$;
                moreover, by construction of $\alpha$, $g\inv[y]$ contains $f\inv[z]$ for some $z\in Z$.
                Hence, $|g\inv [y]|>n$. Finally, if $|f\inv
                [y]|<n$, then there must exist $z\in f\inv[y]$ with
                $1< |(\alpha\circ f)\inv[z]|\leq n$. By construction of $\alpha$,
                $|(\alpha\circ f)\inv[z]|=n$ is impossible, so this can occur at most $s_f(>1,<n)\leq\nu_f$ times. Therefore,
                $s_g(n)\leq\nu_f$. Now Lemma \ref{LEM:UBOUND:FINITE:sg}
                implies that $s_h(n)\leq\nu_f$ for all $h\in\gwithS$; whence, $f\nin\gwithS$,
                contradicting its $\S$-minimality.
            \end{proof}

        \subsection{The decline of the valley of giants}

            \begin{lem}\label{LEM:STRSUPP:s'decreasing}
                If $f$ is $\S$-minimal, then it satisfies (s'dec), i.e.
                $s_f'$ is strictly decreasing.
            \end{lem}
            \begin{proof}
                Assume there exist
                $\psi_1<\psi_2$ in the strong support of $s_f$ with $s_f(\psi_1)\leq
                s_f(\psi_2)$, and let $\psi_1$ be minimal with this property. $\sigma_f=\kappa$ by Lemma
                \ref{LEM:MU:sigmaIsLarge}, and $\rho_f<\kappa$ by Lemma
                \ref{LEM:MU:rhoIsSmall}, so in particular $\mu_f<\psi_1$. Lemma
                \ref{LEM:FINSUPP:nuFinite->sf(n)=0,nuInfinite->sf(n)<=nu}
                tells us that $\psi_1$ cannot be finite.\\
                Fix $Y'_{\psi_2}\subseteq
                Y^f_{\psi_2}$ such that $|Y'_{\psi_2}|=|Y^f_{\psi_1}|$. This is possible since $s_f(\psi_1)\leq
                s_f(\psi_2)$. Choose $S \subseteq Y^f_{\mu_f}$ such that $Y^f_{\mu_f}\sm S$ is
                large and such
                that $|S|=|f\inv [Y'_{\psi_2}]|$.  Let $\alpha$ map
                $Y^f_{\psi_1}\cup S$ bijectively onto $f\inv
                [Y'_{\psi_2}]$, and $Y'_{\psi_2}$ injectively into $f\inv[Y^f_{\psi_1}]$ in
                such a way that for all $y\in Y^f_{\psi_1}$
                there exists $z\in Y'_{\psi_2}$ with $\alpha(z)\in
                f\inv[y]$. We can do that since $|Y'_{\psi_2}|=|Y^f_{\psi_1}|$. Let $\alpha$ moreover map
                all $y\in Y^f_{>\mu_f}\sm (Y^f_{\psi_1}\cup Y'_{\psi_2})$ into $f\inv[y]$. Domain and
                range of $\alpha$ are
                co-large as they are disjoint from $Y^f_{\mu_f}\sm S$ and $f\inv[Y^f_{\mu_f}\sm S]$, respectively,
                so the function can be extended to a bijection on
                $X$; set $g=f\circ\alpha\circ f$.\\
                We claim that $s_g(\psi_1)=0$ and calculate $|g\inv
                [y]|$ for all $y\in X$. If $y\in Y^f_\xi$,
                where $\xi<\psi_1$ is infinite, then $|g\inv
                [y]|\leq\xi<\psi_1$ since $(f\circ\alpha)\inv[y]\sub
                Y^f_{\xi}\cup Y^f_{\mu_f}$ by construction of
                $\alpha$. If $\xi<\psi_1$ is finite, then for the same
                reason we have that $|g\inv[y]|$ is finite, so again
                $|g\inv[y]|<\psi_1$. If $y\in Y^f_{>\psi_1}\sm
                Y'_{\psi_2}$, then $\alpha(y)\in f\inv[y]$ and so
                $|g\inv[y]|\geq |f\inv[y]|>\psi_1$, and if $y\in
                Y'_{\psi_2}$, then $f\inv[y]\cap (\alpha\circ
                f)[X]=f\inv[y]$ by construction of $\alpha$, so again
                $|g\inv[y]|\geq |f\inv [y]|>\psi_1$. Finally, consider
                $y\in Y^f_{\psi_1}$. Then by construction of $\alpha$ there exists
                $z\in Y'_{\psi_2}$ with $\alpha(z)\in f\inv [y]$. But
                then $|g\inv [y]|\geq |f\inv[z]|=\psi_2> \psi_1$, and we have shown $s_g(\psi_1)=0$.\\
                Because $s_f'$ is strictly decreasing below $\psi_1$ by the choice of
                $\psi_1$, its support below $\psi_1$ is finite;
                therefore, unless $\psi_1=\aleph_0$, there exists an infinite $\lambda<\psi_1$ such that
                $\supp'(s_f)\cap (\lambda,\psi_1)$ is empty; moreover, if $\psi_1>\nu_f$, we can certainly
                choose $\lambda$ so that $\lambda\geq\nu_f$.\\
                Consider the case where $\psi_1>\nu_f$ and $\psi_1\neq\aleph_0$; in that
                case, $s_f$ vanishes on the interval
                $(\lambda,\psi_1)$. Because $f$ is $\S$-minimal we must have $\nu_g=\nu_f<\psi_1$. Since $s_g(\xi)\leq
                s_f(>\lambda,\leq\xi)+\min(\nu_f,s_f(>\xi,\leq\nu_f))=s_f(>\lambda,\leq\xi)=0$
                for all $\lambda <\xi< \psi_1$ by Lemma \ref{LEM:UBOUND:singular},
                the same lemma implies that
                $s_h(\psi_1)\leq s_g(>\lambda,\leq\psi_1)+\min(\nu_g,s_g(>\psi_1,\leq\nu_g))=
                s_g(>\lambda,\leq\psi_1)=0$ for all $h\in\gwithS$, so $f\nin\gwithS$, in contradiction
                with the $\S$-minimality of $f$.\\
                If $\psi_1\leq\nu_f$ and $\psi_1\neq\aleph_0$, then $s_f(\psi_1)>\nu_f$ as
                $\psi_1\in\supp'(s_f)$; also, $\nu_g=\nu_f$ by the $\S$-minimality of $f$. By Lemma
                \ref{LEM:UBOUND:singular}, we have $s_g(\xi)\leq s_f(>\lambda, \leq \xi)+\nu_f\leq \nu_f$
                for all $\lambda<\xi<\psi_1$. Therefore by the same lemma, if $h\in\gwithS$, then
                $s_h(\psi_1)\leq s_g(>\lambda, \leq \psi_1)+\nu_g\leq\nu_g$.
                But $\nu_g= \nu_f < s_f(\psi_1)$, so $h\neq f$, and we cannot get $f$ back
                from $g$ and $\S$, again contradicting that $f$ be
                $\S$-minimal.\\
                Finally, if $\psi_1=\aleph_0$, then $s_h(\psi_1)=0<s_f(\psi_1)$
                for all $h\in\gwithS$ by Lemma \ref{LEM:UBOUND:regular}, finishing the last case.
            \end{proof}

        \subsection{The king}

            \begin{lem}\label{LEM:EPSILON:eps=1ORinfinite}
                If $f$ is $\S$-minimal, then it satisfies ($\varepsilon$), i.e. $\ef=1$ or $\ef$ is infinite.
            \end{lem}
            \begin{proof}
                Assume not, and fix $y\in Y^f_{\varepsilon_f}$ and
                any $z\neq y$. Let $\alpha\in\S$ be so that it maps
                $\{y,z\}$ injectively into $f\inv [y]$. Then setting
                $g=f\circ\alpha \circ f$ we have that $|g\inv [y]|\geq |f\inv[y]\cup f\inv[z]|>
                \varepsilon_f$. All functions generated by $g$ with $\S$
                have a class larger than $\varepsilon_f$, which
                implies $f\nin\gwithS$ and contradicts that $f$ is
                $\S$-minimal.
            \end{proof}

            \begin{lem}\label{LEM:EPSILON:BEYONDNU:sf(epsilon)Infinite}
                Let $f$ be $\S$-minimal.
                Then it satisfies (\#$\varepsilon$), i.e.
                if $\nu_f<\varepsilon_f$, then $s_f(\varepsilon_f)$ is infinite.
            \end{lem}
            \begin{proof}
                By Lemma \ref{LEM:STRSUPP:s'decreasing}, the
                restriction of $s_f$ to its support beyond $\nu_f$ is
                strictly decreasing, so the support beyond $\nu_f$ is finite and thus $s_f(\varepsilon_f)>0$.
                Assume $s_f(\varepsilon_f)<\aleph_0$. By Lemma
                \ref{LEM:EPSILON:eps=1ORinfinite}, $\varepsilon_f$ is one or infinite; $\varepsilon_f=1$, however, is
                clearly impossible since it would mean that $f$ is injective but has only finitely many kernel classes.
                Fix $S\subseteq Y^f_{\mu_f}$ such
                that $|S|=|f\inv[Y^f_{\varepsilon_f}]|\geq\ef$ and such that $Y^f_{\mu_f}\sm S$ is
                still large. Let $\alpha$ map $S$ bijectively onto
                $f\inv[Y^f_{\varepsilon_f}]$, as well as $Y^f_{\varepsilon_f}$ injectively into
                $f\inv[S]$.
                The domain of $\alpha$ is disjoint from
                $Y^f_{\mu_f}\sm S$ and hence co-large, and so is its range as
                it is disjoint from $f\inv[Y^f_{\mu_f}\sm S]$, so $\alpha$ can
                be extended to a permutation on $X$. The function $g=f\circ\alpha\circ f$
                satisfies $s_g(\geq\ef)\geq s_f(\varepsilon_f)+1$. Indeed, if $y\in
                Y^f_{\varepsilon_f}$, then $|g\inv [y]|\geq
                |f\inv[y]\cap (\alpha\circ f)[X]|=|f\inv[y]|=\varepsilon_f$
                since $f\inv[y]\cap (\alpha\circ f)[X]=f\inv[y]$ by construction of $\alpha$. Also, taking an arbitrary
                $z\in Y^f_{\varepsilon_f}$ and setting $w=f\circ\alpha(z)\in S$, we have
                $|g\inv[w]|\geq\ef$. Thus indeed, $s_g(\geq\ef)\geq
                |Y^f_{\ef}\cup\{w\}|=s_f(\ef)+1$. However, Lemma
                \ref{LEM:UPPERBOUND:sg(>xi)} gives us $s_g(>\ef)=0$,
                and so $s_g(\varepsilon_f)\geq s_f(\varepsilon_f)+1$.
                Since $\nu_g=\nu_f<\ef$ by the $\S$-minimality of $f$, we have that
                Lemma \ref{LEM:LBOUND:BEYONDNU:sg(xi)lowerbound}
                yields $s_h(\varepsilon_f)\geq s_g(\varepsilon_f) >
                s_f(\varepsilon_f)$ for all $h\in\gwithS\sm\S$, so
                $f\nin\gwithS$, contradicting the assumption that $f$ is $\S$-minimal.
            \end{proof}

            \begin{lem}\label{LEM:EPSILON:epsilonRegularOrMaximum}
                Let $f$ be $\S$-minimal. Then it satisfies ($\varepsilon$reg),
                i.e., either
                $s_f(\varepsilon_f)>0$ or $\ef$ is regular.
            \end{lem}
            \begin{proof}
                Assume $\ef$ is singular and that $s_f(\ef)=0$.
                Let $\eta<\varepsilon_f$ be the cofinality of
                $\varepsilon_f$, let $\vartheta\geq\eta$ be in the
                support of $f$, and fix $y\in Y^f_\vartheta$. Let
                $(\zeta_\tau)_{\tau<\eta}$ be a strictly increasing sequence of
                cardinalities in the support of $s_f$ which is
                cofinal in $\varepsilon_f$ and larger than $\mu_f$, and fix
                $y_\tau\in Y^f_{\zeta_{\tau}}$ for all $\tau<\eta$. Set
                $Y=\{y_\tau: \tau<\eta\}$. Let $\alpha$ map $Y$ injectively into
                $f\inv[y]$, and extend it to a bijection. This is possible
                since $\alpha$ is not defined on $Y^f_{\mu_f}$ and its range is disjoint
                from $f\inv [Y^f_{\mu_f}]$ and since the two sets are
                large by Lemma \ref{LEM:MU:sigmaIsLarge}. Set $g=f\circ \alpha\circ
                f$. Then $|g\inv[y]|\geq
                |\bigcup_{\tau<\eta}f\inv[y_\tau]|=\sum_{\tau<\eta}\zeta_\tau=\varepsilon_f$.
                Therefore, $g$ has a kernel class larger than all
                kernel classes of $f$, so that it cannot generate $f$
                together with $\S$, contradiction.
            \end{proof}

        \subsection{Continuity}

            \begin{lem}\label{LEM:MIDSUPP:noJumpsAtSingular}
                Let $f$ be $\S$-minimal. Then
                $s_f(\geq\xi)=\min\{s_f(\geq\zeta):\zeta<\xi\}$ for
                all singular $\xi\leq\chi_f$.
            \end{lem}
            \begin{proof}
                Assume there is $\xi\leq\chi_f$ singular with
                $s_f(\geq\xi)<\vartheta=\min\{s_f(\geq\zeta):\zeta<\xi\}$, and
                let $\eta<\xi$ be the cofinality of $\xi$. Clearly,
                $\xi>\mu_f$. Let $\max\{\eta,\mu_f\}<\zeta<\xi$ be so that
                $s_f(\geq\zeta)=\vartheta$.\\
                Observe next that for all
                $\zeta\leq\psi<\xi$ and all $\lambda<\vartheta$ there
                exists $\psi\leq \psi'<\xi$ with
                $s_f(\psi')>\lambda$, for otherwise $\vartheta=s_f(\geq
                \psi)\leq \lambda\cdot\xi$, implying $\vartheta=\xi$,
                and thus $\ef\leq\vartheta=\xi\leq\ef$. However,
                $\vartheta=\ef$ implies $s_f(\ef)=0$ since $\zeta<\chi_f$ and $s_f(\geq\zeta)=\vartheta$,
                contradicting Lemma
                \ref{LEM:EPSILON:epsilonRegularOrMaximum}. By our observation we
                can thin out the interval $(\zeta,\xi)$ and
                find
                a strictly increasing sequence of
                cardinals $(\zeta_\tau)_{\tau<\eta}$ greater than $\zeta$ and cofinal in
                $\xi$, such that that the sequence
                $(\delta_\tau)_{\tau<\eta}=(s_f(\zeta_\tau))_{\tau<\eta}$
                is increasing and has the property that for all $\lambda<\vartheta$ there exists $\tau<\eta$ such that
                $\delta_\tau>\lambda$. Write
                $Y^f_{\zeta_\tau}=\{y_{\zeta_\tau}^i:i<\delta_\tau\}$ for all $\tau<\eta$ (with the variable $i$
                referring to all ordinals below $\vartheta$).
                Set $S^i=\{y_{\zeta_\tau}^i: \tau<\eta \wedge
                i<\delta_\tau \}$, for all $i<\vartheta$.\\
                Fix a set $Z\sub Y^f_{\geq\zeta,<\xi}$ such that
                $|Z|=\vartheta$ and write $Z=\{z_i:i<\vartheta\}$ (again with $i$ referring to ordinals).
                Let $\alpha$ map $S^i$ injectively
                into $f\inv[z_i]$, for all $i<\vartheta$,
                extend $\alpha$ to a bijection, and set
                $g=f\circ\alpha\circ f$. Then $|g\inv
                [z_i]|\geq
                |f\inv[S^i]|=|\bigcup \{f\inv [y_{\zeta_\tau}^i]: \tau<\eta \wedge
                i<\delta_\tau \}|=\sum_{\tau<\eta\wedge i<\delta_\tau}\zeta_\tau=\xi$, the latter
                equality holding since the condition $i<\delta_\tau$ only cuts away an initial segment of the
                sequence $(\zeta_\tau)_{\tau<\eta}$ which is cofinal in $\xi$.
                Thus, $s_g(\geq\xi)\geq \vartheta$. Now $g$ does not have
                any kernel class larger than all kernel classes of $f$,
                because $f$ is $\S$-minimal; hence, $s_g(\geq\xi)$ is larger that all cardinals
                in $\supp(s_g)$, and thus $\xi<\chi_g$. Moreover,
                $g$ satisfies ($\varepsilon$) and ($\varepsilon$reg),
                by Lemmas \ref{LEM:EPSILON:eps=1ORinfinite} and
                \ref{LEM:EPSILON:epsilonRegularOrMaximum}. Therefore,
                $s_h(\geq\xi)\geq s_g(\geq\xi)=\vartheta>s_f(\geq\xi)$ for all
                $h\in\gwithS\sm\S$ by Lemma
                \ref{LEM:MIDSUPP:LBOUND:sf(>=xi)}, contradicting
                that $f$ is $\S$-minimal.
            \end{proof}
            \begin{lem}\label{LEM:MIDSUPP:noJumpsAtFinite}
                Let $f$ be $\S$-minimal. Then $s_f(\geq n)=s_f(\geq 2)$ for all
                finite $2\leq n\leq\chi_f$.
            \end{lem}
            \begin{proof}
                It suffices to show that $s_f(n)=s_f(n+1)$ for all finite $2\leq n<\chi_f$.
                Assume to the contrary that $s_f(\geq n)>s_f(\geq n+1)$ for some finite $2 \leq n<\chi_f$.
                By Lemma \ref{LEM:EPSILON:eps=1ORinfinite}, $f$
                satisfies ($\varepsilon$). This, together with the fact that there
                is no $\zeta\in\supp(s_f)$ with $\zeta\geq s_f(\geq
                n)$, implies that
                $s_f(\geq n)$ must be infinite, and hence $s_f(\geq n)=s_f(n)$ as
                $s_f(>n)<s_f(\geq n)$. Let $\alpha$ map $Y^f_{n}$ injectively into $f\inv[Y^f_{n}]$ in
                such a way that $|f\inv [y]\cap \alpha[Y^f_{n}]|=2$ for all $y\in Y^f_{n}$.
                Because $\alpha$ satisfies ($\mu$) and ($\sigma$) by Lemmas \ref{LEM:MU:muIsOneOrInfinite} and
                \ref{LEM:MU:sigmaIsLarge}, we have that $Y_1^f\neq Y_n^f$ is large, so we
                can extend $\alpha$ to a permutation of $X$ and set $g= f\circ \alpha\circ f$.
                Then for all $y\in
                Y^f_n$ we have that $|g\inv[y]|\geq 2\cdot n>n$.
                Hence, $s_g(\geq n+1)\geq s_f(\geq n)> s_f(\geq n+1)$.
                We clearly have $\varepsilon_g=\ef$ and
                $s_g(\ef)=0$ iff $s_f(\ef)=0$, as $f$ is $\S$-minimal, so
                $s_g(\geq n+1)\geq s_f(\geq n)$ and $\chi_f>n$ imply
                $\chi_g>n+1$. Also, $g$ satisfies ($\varepsilon$reg) by Lemma
                \ref{LEM:EPSILON:epsilonRegularOrMaximum}. Thus,
                Lemma \ref{LEM:MIDSUPP:LBOUND:sf(>=xi)} implies
                that $s_h(\geq n+1)\geq s_g(\geq n+1)> s_f(\geq n+1)$ for all
                $h\in\gwithS\sm\S$, so $f$ is not $\S$-minimal, contradiction.
            \end{proof}

            \begin{lem}\label{LEM:SCONT}
                Let $f$ be $\S$-minimal. Then it satisfies \cscont,
                i.e. $s_f(\geq\xi)=\min\{s_f(\geq\zeta):\zeta<\xi\}$ for all singular
                $\xi\leq\chi_f$ and $s_f(\geq n)=s_f(\geq 2)$ for all finite $2\leq n\leq\chi_f$.
            \end{lem}
            \begin{proof}
                This is the consequence of Lemmas
                \ref{LEM:MIDSUPP:noJumpsAtSingular} and
                \ref{LEM:MIDSUPP:noJumpsAtFinite}.
            \end{proof}

        \subsection{The rage of the lion-tail}

            \begin{lem}\label{LEM:WSUPP:nu=kappa->sf(kappa)=0}
                If $f$ is $\S$-minimal and nonconstant, then it satisfies ($\kappa$), i.e.,
                if $\nu_f=\kappa$, then $s_f(\kappa)=0$.
            \end{lem}
            \begin{proof}
                Assume $\nu_f=\kappa$ and $s_f(\kappa)>0$ and let $y\in Y^f_\kappa$.
                Let $\alpha$ map $f[X]$ injectively
                into a co-large part of $f\inv[y]$. Since both domain
                and range of $\alpha$ are co-large, we can extend it
                to a function in $\S$. Then $g=f\circ\alpha\circ f$ is
                constant and generates together with $\S$ a proper subclone of
                $\cl{\{f\}\cup\S}$, contradicting that $f$ is $\S$-minimal.
            \end{proof}

            \begin{lem}\label{LEM:CHI:sf(>=chi)Finite}
                If $f$ is $\S$-minimal, then it satisfies ($\chi$), i.e. if $\ef\leq\nu_f$,
                then $s_f(\geq\chi_f)$ is finite.
            \end{lem}
            \begin{proof}
                We can assume that $f$ is nonconstant and that $\varepsilon_f>1$; then $\varepsilon_f$ is
                infinite by Lemma \ref{LEM:EPSILON:eps=1ORinfinite}. Also, we may assume that $\mu_f<\chi_f$,
                for otherwise $s_f(\kappa)>0$ as $s_f(\mu_f)=\sigma_f=\kappa$ by Lemma \ref{LEM:MU:sigmaIsLarge} and thus
                $\nu_f\geq\ef=\kappa$,
                contradicting that $f$ satisfies ($\kappa$) by Lemma \ref{LEM:WSUPP:nu=kappa->sf(kappa)=0}.
                Suppose $s_f(\geq\chi_f)$ is infinite;
                we want to derive a contradiction. By Lemma \ref{LEM:SCONT}, $f$ satisfies
                (scont), and therefore $\chi_f\leq 2$ or $\chi_f$ is infinite and
                regular by Lemma \ref{LEM:CHI:chiIsRegular}.
                Fix $y\in X$ with $|f\inv[y]|\geq s_f(\geq\chi_f)$.
                Let $\alpha$ map
                $Y^f_{\geq\chi_f}$ injectively into $f\inv[y]$, and
                a suitable part of $X\sm f[X]$ bijectively onto
                $f\inv[Y^f_{\geq\chi_f}\sm\{y\}]$; this is possible as $|f\inv[Y^f_{\geq\chi_f}]|\leq
                \ef\cdot |Y^f_{\geq\chi_f}|\leq\ef\leq\nu_f$. Extend
                $\alpha$ to a bijection on $X$. The function
                $g=f\circ\alpha\circ f$ satisfies
                $s_g(\geq\chi_f)\leq 1$. Indeed, if
                $|g\inv[z]|\geq\chi_f$, then either
                $|f\inv[z]|\geq\chi_f$ or there exists $w\in f\inv[z]$
                with $|(\alpha\circ f)\inv [w]|\geq\chi_f$, because
                $\chi_f\leq 2$ or $\chi_f$
                is infinite and regular. But if
                $|f\inv[z]|\geq\chi_f$, then for $z\neq y$ we have that $f\inv[z]\cap \alpha\circ
                f[X]=\emptyset$, by definition of $\alpha$, so $|g\inv[z]|=0$; the
                other possibility
                does not occur unless $z=y$, and we have shown $s_g(\geq\chi_f)\leq 1$. Therefore,
                $s_h(\geq\chi_f)$ is finite for all $h\in\gwithS$ by Lemma \ref{LEM:UBOUND:sg(>=xi)},
                and hence $f\nin\gwithS$, contradiction.
            \end{proof}

        \subsection{Existence of the hole}

            \begin{lem}\label{LEM:WSUPP:removalOfBoundedElements}
                Let $f\in\Oo$ satisfy ($\mu$), ($\nu$), ($\sigma$), (s'dec), and (n).
                Then there exists $g\in\fwithS$ such that $s_g(\xi)=0$ for all $\xi<\efp$ with $\xi\nin\supp'(s_f)$,
                and $s_g(\xi)=s_f(\xi)$ for all other $\xi\leq \kappa$. In particular, there are no
                elements below $\varepsilon_f'$ in the weak support of
                $s_g$.
            \end{lem}
            \begin{proof}
                This can be proven exactly like Lemma
                \ref{LEM:WSUPP:removalOfBoundedElementsLambda}, replacing
                $\lambda_f'$ by $\efp$.
            \end{proof}

            \begin{lem}\label{LEM:NEC:(lambda')}
                Let $f$ be $\S$-minimal. Then $f$ satisfies ($\lambda$'), i.e. if $\ef>\nu_f$, then $s_f(\xi)=0$ for all
                $\xi\in (\lambda_f',\nu_f]$.
            \end{lem}
            \begin{proof}
                Assume $\ef=\efp>\nu_f$. By Lemmas \ref{LEM:MU:muIsOneOrInfinite},
                \ref{LEM:NU:(mu=1&NuInfinite)ORMuInfinite},
                \ref{LEM:MU:sigmaIsLarge},
                \ref{LEM:STRSUPP:s'decreasing} and
                \ref{LEM:FINSUPP:nuFinite->sf(n)=0,nuInfinite->sf(n)<=nu},
                $f$ satisfies the conditions of Lemma
                \ref{LEM:WSUPP:removalOfBoundedElements}. Therefore
                there exists $g\in \fwithS\sm\S$ such that $s_g(\xi)=0$
                for all $\xi\in (\lambda_f',\nu_f]$. By Lemma \ref{LEM:UBOUND:singular},
                $s_h(\xi)\leq s_g(\xi)+\min(\nu_g,s_g(>\xi,\leq\nu_g))=0+\min(\nu_g,0)=0$
                for all $\xi\in (\lambda_f',\nu_f]$ and
                all $h\in\gwithS$, so in particular this holds for $f$.
            \end{proof}

    \section{Proofs of the corollaries}

        \begin{proof}[Proof of Theorem \ref{THM:whenEqual}]
            Assume first that $\fwithS=\gwithS$. By Lemma
            \ref{LEM:STABLE:MU} we have $\mu_g=\mu_f$,
            $\nu_g=\nu_f$, and $\varepsilon_g=\ef$.
            Lemma \ref{LEM:STRSUPP:s'dec->stable} implies
            $s_g'=s_f'$. We have $\chi_g=\chi_f$ by Lemma
            \ref{LEM:CHI:ConditionsForChiStable}, and by Lemma
            \ref{LEM:MIDSUPP:sf(>=xi)StableConditions} we have
            $s_g(\geq\xi)=s_f(\geq\xi)$ for all $\xi<\chi_f$. Obviously $s_f(\ef)>0$ implies $s_g(\ef)>0$.\\
            For the other direction, assume first that
            $\ef>\nu_f$. By Lemma \ref{LEM:WSUPP:BELOWLAMBDA:addingOfClasses},
            there exists $h\in\fwithS$ such that
            $s_h\rest_{[1,\lambda_f')}=s_g\rest_{[1,\lambda_f')}$
            and $s_h\rest_{[\lambda_f',\kappa]}=s_f\rest_{[\lambda_f',\kappa]}$;
            thus, $s_h\rest_{(\nu_f,\kappa]}=s_g\rest_{(\nu_f,\kappa]}$ as $s_g'=s_f'$.
            Also, we have that
            $\supp(s_h)\cap(\lambda_f',\nu_f]=\supp(s_f)\cap(\lambda_f',\nu_f]$ is empty, by \clamp;
            for the same reason, $s_g$ vanishes in that interval, too. Therefore,
            $s_h=s_g$ so that since by Lemma \ref{LEM:NU:(nu)->NuStable} also $\nu_h=\nu_f=\nu_g$ we
            conclude $g\in\fwithS$.\\
            Next assume $\ef\leq\nu_f$ and $s_f(\geq \chi_f)=0$; then also $s_g(\geq \chi_f)=0$ as
            $\varepsilon_g=\ef$ and since $s_g(\varepsilon_g)=0$ iff $s_f(\ef)=0$.
            By Lemma \ref{LEM:MIDSUPP:MODIFY:arbitraryBelowChi} there exists $h\in\fwithS$ such that
            $s_h\rest_{[1,\chi_f)}=s_g\rest_{[1,\chi_f)}$, and such that
            $s_h(\xi)=s_f(\xi)=0$ for all $\xi\geq\chi_f$; hence,
            $s_h=s_g$ and we are done.\\
            Finally, if $\ef\leq\nu_f$ and $s_f(\geq\chi_f)>0$, then $s_f(\geq\chi_f)$ is
            finite by \cchi, and so is $s_g(\geq\chi_g)$ for the same reason. Also, $s_g(\geq\chi_g)>0$ as
            $\varepsilon_g=\ef$ and since $s_g(\varepsilon_g)=0$ iff $s_f(\ef)=0$. With
            the help of Lemma \ref{LEM:TAIL:MODIFY:arbitraryEveryWhere}
            we can construct $h\in\cl{\{f\}\cup\S}$
            such that $s_h=s_g$.
        \end{proof}

        \begin{proof}[Proof of Corollary \ref{COR:numberOfMinimal}]
            By Theorem \ref{THM:whenEqual}, the clone an
            $\S$-minimal function $f$ generates is fully determined by
            the decreasing sequences $s_f'(\xi)$ and $s_f(\geq\xi)$, as well as by the values $\mu_f$,
            $\nu_f$, $\chi_f$, $\ef$, and $s_f(\ef)$. Since $s_f'(\xi)$ and $s_f(\geq\xi)$
            are decreasing, they are determined by the finitely many
            points where they decrease, together with their values at
            those points. Therefore, for all determining
            parameters we have at most as many possibilities as
            there are cardinals below $\kappa=\aleph_\alpha$,
            which is $\max\{|\alpha|,\aleph_0\}$, so the number of
            clones minimal in $[\cl{\S},\O]$ is not more than that.\\
            On the other hand, using Theorem \ref{THM:mainTheorem}
            one sees that the functions $f\in\Oo$ with
            $\mu_f=\kappa$, $s_f(\kappa)=\kappa$ and
            $\nu_f=\nu<\kappa$ are $\S$-minimal for all
            $\nu<\kappa$, and by Theorem \ref{THM:whenEqual} they
            generate distinct clones. Therefore, the number of
            clones minimal in $[\cl{\S},\O]$ is at least
            $\max\{|\alpha|,\aleph_0\}$.
        \end{proof}

        \begin{proof}[Proof of Corollary \ref{COR:countableCase}]
            The $\S$-minimality of the functions which generate those monoids can easily be verified by Theorem
            \ref{THM:mainTheorem}.\\
            To see that the mentioned monoids are the only monoids
            minimal in $[\S,\Oo]$,
            let $f$ be $\S$-minimal and non-constant. If $\mu_f=\aleph_0$ and
            $\nu_f <\aleph_0$, then $f$ with $\S$ generates
            $\I_{\nu_f}$. We cannot have $\mu_f=\aleph_0$ and $\nu_f=\aleph_0$, because this would contradict
            \cchi or \csig. So let $\mu_f=1$; then $\nu_f$ is zero or infinite by ($\nu$). We
            distinguish two cases. Assume first that $\ef=\mu_f=1$.
            Then $\nu_f >0$ since $f\nin\S$, so $\nu_f$ is infinite and it is easily seen that
            in this case, $f$ generates $\H$. Now consider the
            case where $\ef>1$; we claim that this cannot happen. Indeed, we would have to have
            $\ef=\aleph_0$ by ($\varepsilon$).
            By ($\rho$), $s_f(>1)$ is finite and therefore $s_f(\ef)>0$.
            But then $\chi_f=1$ by definition, contradicting \cchi or \csig.
        \end{proof}

\end{document}